\documentclass[12pt]{amsart}
\pdfoutput=1 
\usepackage{microtype}
\usepackage{amsmath, amsthm, graphicx, amssymb,verbatim,pst-all,color}
\usepackage[margin=1.0in,footskip=0.45in]{geometry}
\usepackage{mathptmx}
\usepackage{multirow}
\usepackage{subfig}
\usepackage{enumerate}
\usepackage{array,nicefrac}
\usepackage{url}
\usepackage{times}
\usepackage{textcomp}
\usepackage{verbatim}
 \usepackage[T1]{fontenc}
\usepackage{epsf}
\usepackage{stmaryrd}
\usepackage{cases} 
\usepackage{relsize}
\usepackage{mathtools}
\usepackage{amscd, amsfonts}
\usepackage[shortlabels]{enumitem}
\usepackage[all]{xy}
\usepackage{lscape}
\usepackage{tikz-cd}
\usepackage{soul}
\usepackage{pgf}
\usepackage[hidelinks]{hyperref}
\ifx\JPicScale\undefined\def\JPicScale{1.0}\fi
\unitlength \JPicScale mm
{
\theoremstyle{plain}
  \newtheorem{theorem}{Theorem}[section]
  \newtheorem{proposition}[theorem]{Proposition}   
  \newtheorem{lemma}[theorem]{Lemma}

  \newtheorem{corollary}[theorem]{Corollary}
}
{
\theoremstyle{definition}

  \newtheorem{example}[theorem]{Example}
  \newtheorem{definition}[theorem]{Definition}
  \newtheorem{remark}[theorem]{Remark}
  \newtheorem{conv}[theorem]{Convention}
  
}
\newcommand{\Gm}{\mathbb{C}^{*}}

\newcommand{\Bl}{\operatorname{Bl}}

\newcommand{\cL}{\mathcal{L}}
\newcommand{\cO}{\mathcal{O}}
\newcommand{\Proj}{\operatorname{Proj}}
\newcommand{\Sym}{\operatorname{Sym}}

\usepackage[symbol]{footmisc}

\newcommand{\CC}{\mathbb{C}}

\newcommand{\ZZ}{\mathbb{Z}}
\newcommand{\NN}{\mathbb{N}}

\newcommand{\Aff}{\mathbb{A}}
\newcommand{\PP}{\mathbb{P}}
\newcommand{\SL}{\text{SL}}

\newcommand{\quotient}{ / \! \! /}

\newcommand{\wbuildtdn}{\mathcal{H}_{\mathcal{A}}}

\newcommand{\wtdn}{T_{d,n}^{\mathcal{A}}}

\newcommand{\wtdno}{(T_{d,n}^{\mathcal{A}})^{\circ}}

\newcommand{\bP}{\mathbb{P}}

\author{Patricio Gallardo, Javier Gonz\'alez Anaya, Jos\'e Luis Gonz\'alez \and Evangelos Routis}
\address{
{\small Department of Mathematics,
University of California, Riverside,
900 University Ave.,
Riverside, CA 92521,
United States.
}
}
\email{pgallard@ucr.edu}

\address{
{\small Department of Mathematics and Computer Science, 
Santa Clara University, 
500 El Camino Real,
Santa Clara, CA 95053,
United States.}
}
\email{jgonzalezanaya@scu.edu}

\address{
{\small Department of Mathematics,
University of California, Riverside,
900 University Ave.,
Riverside, CA 92521, 
United States.}
}
\email{jose.gonzalez@ucr.edu}

\address{
{\small 
Mathematics Institute,
Zeeman Building,
University of Warwick,
Coventry CV4 7AL, United Kingdom.}
}
\email{evangelos\_routis@alumni.brown.edu}

 \let\oldbibliography\thebibliography
\renewcommand{\thebibliography}[1]{%
  \oldbibliography{#1}%
  \setlength{\itemsep}{4pt}%
  \setlength{\lineskiplimit}{-\maxdimen}
}

\begin{document}
\title{
Higher-dimensional Losev-Manin spaces and their geometry
}

\begin{abstract}
The classical Losev-Manin space is a toric compactification of the moduli space of $n$ points in the affine line modulo translation and scaling. 
Motivated by this, we study its higher-dimensional toric counterparts, which compactify the moduli space of $n$ distinct labeled points in affine space modulo translation and scaling. We show that these moduli spaces are a fibration over a product of projective spaces---with fibers isomorphic to the Losev-Manin space---and that they are isomorphic to the normalization of a Chow quotient. Moreover, we present a criterion to decide whether the blow-up of a toric variety along the closure of a subtorus is a Mori dream space. 
As an application, we demonstrate that a related generalization of the moduli space of pointed rational curves constructed by Chen, Gibney, and Krashen is not a Mori dream space when the number of points is at least nine, regardless of the dimension.
\end{abstract}

\maketitle


\section{Introduction}
Toric compactifications of moduli spaces are fundamental objects in algebraic geometry, having deep connections to tropical geometry \cite{tevelev2007compactifications}, 
logarithmic geometry \cite{chen2013chow, ascher2016logarithmic}, 
and combinatorial constructions such as secondary and fiber polytopes
\cite{craw2007fiber, kapranov1991quotients}. Several well-known toric moduli spaces include Alexeev's work on stable pairs \cite{alexeev2002complete} and toric Hilbert schemes \cite{peeva2002toric, chuvashova2008main}. 
This article is concerned with a recent application of toric compactifications which aims to compactify moduli spaces of $n$ distinct labeled points in toric varieties of arbitrary dimension over the complex numbers. 
Examples using related techniques to ours include recent work by Schaffler-Tevelev~\cite{schaffler2022compactifications} and Schaffler-Di~Rocco~\cite{di2022families} on the moduli space of points in toric varieties up to the action of the maximal torus.  We remark that other common techniques used to compactify moduli spaces, like the MMP and K-stability, do not apply to our case, as these methods are geared towards compactifying moduli spaces of varieties possibly marked by divisors, and in our case the varieties are marked by higher codimensional loci. 
Throughout this article, we follow the conventions of \cite{Fulton_Toric_Varieties_Book} for the definitions and notation related to toric varieties.
In particular, the toric varieties in this article are normal by assumption. 

A central example is the Losev-Manin compactification $\overline{M}_{0,n}^{LM}$ \cite{losev2000new-moduli-spaces}. This toric variety, associated with the permutahedron, has two modular interpretations: As a compactification of the moduli space of $n$ points in $\mathbb{P}^1$ modulo projective linear transformations, or as a compactification of the moduli space of $(n-1)$ points in the affine line modulo translation and scaling. The second interpretation leads us to consider the moduli space of $n$ distinct labeled points in affine space $\mathbb{A}^d$ up to translation and scaling.  This last moduli problem was initially studied by Chen, Gibney, and Krashen, and it has a non-toric compactification denoted as $T_{d,n}$, which serves as a generalization of $\overline{M}_{0,n}$ since $T_{1,n} \cong \overline{M}_{0,n+1}$; see \cite{Chen-Gibney-Krashen}.
In recent work \cite{GR}, the authors extended this construction by allowing weights on the points, similar to Hassett's compactification of weighted pointed rational curves \cite{Hassett-weighted}. This generalization gives rise to a toric compactification called $T_{d,n}^{LM}$, which is isomorphic to the Losev-Manin compactification when the dimension is $d=1$. For $d > 1$, it becomes a higher-dimensional generalization of the Losev-Manin space. Our main results provide a description of the geometry of $T_{d,n}^{LM}$ and its applications in studying $T_{d,n}$.

\begin{theorem}\label{thm:MainLM}
The smooth toric compactification $T_{d,n}^{ LM}$ of the moduli space of $n$ points in the affine space $\mathbb{A}^d$ up to translation and scaling satisfies the following:
\begin{enumerate}[label=(\roman*),ref=\thetheorem.(\roman*)]
\item\label{Main thm: I} $T_{d,n}^{LM}$ is isomorphic to the normalization of the Chow quotient $(\mathbb{P}^{d})^{n-1} / \! \! /_{Ch} \mathbb{C}^*$ introduced in Definition~\ref{chow: definition chow quotient}.

\item\label{Main thm: II} There is an irreducible closed locus $\delta \subset T_{d,n}^{LM}$ with $\dim(\delta) =d-1$ such that the canonical map $T_{d,n}\to T_{d,n}^{LM}$ factors as 
\begin{align*}
T_{d,n} \longrightarrow \operatorname{Bl}_{\delta} T_{d,n}^{LM}
\longrightarrow T_{d,n}^{LM},
\end{align*}
and $\operatorname{Bl}_{\delta} T_{d,n}^{LM}$ is not a Mori dream space for $n \geq 9$.

\item\label{Main thm: III} $T_{d,n}^{LM}$ is a Zariski locally trivial fibration over $(\PP^{d-1})^{n-1}$, with fiber isomorphic to the Losev-Manin space $\overline{M}_{0,n+1}^{LM}$.
\end{enumerate}
\end{theorem}

The proofs of Theorem~\ref{thm:MainLM} parts (i), (ii) and (iii) are given in Sections~\ref{sec:AsToricChowQuotient}, \ref{sec:NotMDS} and \ref{sec:Fibrations}, respectively. Preliminaries about $T_{d,n}$ and $T_{d,n}^{LM}$ are presented in Section~\ref{sec:Preliminarytdn}. In the rest of the introduction, we will describe the key aspects of our proof, its context, as well as applications.

We construct our toric Chow quotient as a subscheme in a Chow variety, which represents a functor of certain families of cycles over seminormal bases in characteristic zero \cite[Theorem 3.21]{Kollar-chow}. In particular, the Chow variety has a universal family. Therefore, the isomorphism between
$T_{d,n}^{LM}$ and the normalization of our toric Chow quotient improves over the initial constructions in \cite{GR} and \cite{Chen-Gibney-Krashen}. Indeed, it expands the moduli interpretation of $T_{d,n}^{LM}$ since we can pullback the universal family over the Chow variety to obtain a new family over $T_{d,n}^{LM}$. The existence of an algebraic stack that parameterizes families of points and is represented by the Fulton-MacPherson compactification $\mathbb{P}^d[n]$ of the configuration space of $n$ distinct labeled points in $\mathbb{P}^d$ is still conjectural; see \cite{lundkvist2009moduli} for progress on this issue.

In addition, our construction generalizes similar results for the moduli of points on the line \cite[Lemma~5.1]{dolgachev2015configuration}, points in projective space \cite[Section~8]{schaffler2022compactifications}, and in toric varieties \cite{di2022families}. In certain cases, the Chow quotient is isomorphic to the inverse limit of GIT quotients \cite[Proposition~2.4]{baker2015chow}, and carries rich combinatorial information, for example, in relation to the theory of secondary polytopes \cite{kapranov1991quotients}, and tropical geometry \cite{gibney2011equations}. Additionally, the Chow quotient has topological interpretations \cite{hu2005topological, giansiracusa2022chow}, as well as interpretations in terms of stable maps \cite{chen2013chow}.

Next we discuss Theorem~\ref{Main thm: II}, which we leverage to study the finite generation of the Cox ring of the Chen-Gibney-Krashen compactification $T_{d,n}$. 
Cox rings are a generalization of the homogeneous coordinates of toric varieties from \cite{Cox-ring} to a wider class of algebraic varieties, 
including normal complete varieties with a finitely generated divisor class group; for a survey, see \cite{gonzalez2022finite} and references within.
The Cox ring carries substantial information about the geometry of a variety. As an example, a variety is a toric variety if and only if its Cox ring is isomorphic to a polynomial ring \cite[Corollary~2.10]{Hu-Keel}. A normal $\mathbb{Q}$-factorial projective variety with a finitely generated Cox ring is called a Mori dream space (MDS), as this property guarantees that Mori's program can be carried out for any divisor on $X$ \cite[Proposition~1.11]{Hu-Keel}. 

Studying the MDS property for blow-ups $\operatorname{Bl}_e Y$ of toric varieties $Y$ at a general point $e$ has received significant attention, as the MDS property for other important varieties can be reduced to this case; see for example \cite{Castravet-Tevelev, gonzalez2016some, hausen2018blowing, he2019mori,GGK1,GGK2,GGK3,MRSS, Geography, laface2023intrinsic, nonexist, CastravetLafaceTevelevUgaglia2020}. 
Notable examples include the proof  
that $\overline{M}_{0,n}$ is not a Mori dream space for sufficiently large $n$ \cite{Castravet-Tevelev}, along with subsequent refinements \cite{gonzalez2016some,hausen2018blowing}.
Another example is in the proof that the Fulton-MacPherson compactification $\mathbb{P}^d[n]$ is not a Mori dream space for sufficiently large $n$ and $d \geq 1$ \cite{Fulton-MacPherson2022}.
However, there have been limited results for blow-ups of toric varieties along higher-dimensional centers. 
To address this, we rely on the work of \cite{hausen2010cox} to establish the following result (proved in Theorem~\ref{theorem.blow-up.torus.along.subtorus}), which serves as a key tool in such cases.  

\begin{theorem} \label{introduction: blowup theorem quotient} 
Let $X$ be a complete toric variety and let $T$ be a subtorus of the torus $T_X$ of $X$. 
Let $N \supseteq N_T$ be the lattices of one-parameter subgroups of $T_X$ and $T$.
Let $Y$ be a complete toric variety given by a fan in $N/N_T \otimes \mathbb{R}$ whose set of rays are the images of the rays of X under the natural projection $N \otimes \mathbb{R} \rightarrow N/N_T \otimes \mathbb{R}$ and let $e$ be a point in the torus of $Y$. 
Then, 
\[
\operatorname{Cox}(\operatorname{Bl}_{\overline{T}} X) 
\textnormal{ is finitely generated} 
\ \ \ \Longleftrightarrow \ \ \
\operatorname{Cox}\left(\operatorname{Bl}_e Y\right) 
\textnormal{ is finitely generated}. 
\]
\end{theorem}

To apply the above result to our specific case, we begin by noting that there exists a sequence of iterative smooth blow-ups given by:
\[
    T_{d,n} \longrightarrow \operatorname{Bl}_{\delta} T_{d,n}^{LM} \longrightarrow T_{d,n}^{LM} \longrightarrow \mathbb{P}^{nd-d-1},
\]
where $\dim(\delta) = d-1$. Applying the aforementioned theorem to the intermediate step $\operatorname{Bl}_{\delta} T_{d,n}^{LM}$, we obtain the blow-up $\operatorname{Bl}_{e} Y$, where $Y$ is a priori an unknown toric variety and $e$ a point in its torus. We show that $Y$ is a toric GIT quotient of the Fulton-McPherson compactification of $n$ points in $\mathbb{P}^{d}$, denoted as $\overline{P}^{LM}_{d,n+d}$ and introduced by two of the authors in \cite{GR}. Moreover, in \cite[Lemma 3.4]{Fulton-MacPherson2022}, it was subsequently shown by three of the authors that the Cox ring of the blow-up $\operatorname{Bl}_{e} \overline{P}^{LM}_{d,n+d}$ is not finitely generated for $n \geq 9$, and from this Theorem~\ref{Main thm: II} follows.   
In addition, we show that $T_{d,n}$ is a Mori dream space for $n \leq 3$, which together with Theorem~\ref{Main thm: II} gives us the following.
\begin{theorem} 
\label{theorem.Tdn.not.MDS}
$T_{d,n}$ is a Mori dream space for $n \leq 3$ and 
not a Mori dream space for $n \geq 9$.
\end{theorem}
We remark that $T_{d,2}$ is trivially a MDS because it is isomorphic to a projective space, $T_{1,n}$ is a MDS for $n \leq 5$ because it is isomorphic to $\overline{M}_{0,n+1}$, and $T_{2,3}$ is the blow-up of three disjoint lines in $\mathbb{P}^3$, so it is a MDS by \cite[Theorem 4.4]{dumitrescu2017cones}. Therefore, this result addresses the natural question of finding necessary and sufficient conditions on the dimension $d$ and number of points $n$ for $T_{d,n}$ to be a MDS. 

Finally, we turn to Theorem~\ref{Main thm: III}. Here we utilize the fact that $T_{d,n}^{LM}$ can be obtained through an iterative sequence of smooth blow-ups starting from $\mathbb{P}^{d(n-1)-1}$.
In Section~\ref{sec:Fibrations} we start from a concrete rational map from $\mathbb{P}^{d(n-1)-1}$ to $\left(\mathbb{P}^{d-1}\right)^{n-1}$. Then, the result is derived by demonstrating that this map can be resolved by rearranging the iterative blow-up construction $T_{d,n}^{LM} \longrightarrow \mathbb{P}^{d(n-1)-1}$, while keeping track of the blow-up process. 
This structural result can be used to study the geometry of the higher-dimensional Losev-Manin spaces $T_{d,n}^{LM}$. For example, in Corollary~\ref{cor:DC} we invoke \cite[Theorem~1.3]{costa2011derived} and apply this result to show that if $d \geq 1$ and $n \leq 4$, then $T_{d,n}^{LM}$ has a full strong exceptional collection of line bundles.

\begin{conv}
    Throughout this article we work over the field of complex numbers $\CC$.
\end{conv}
    
\subsection{Acknowledgements}
Patricio Gallardo, Javier Gonz\'alez Anaya and Jos\'e Gonz\'alez are grateful for the working environment at the University of California, Riverside, where part of this research was conducted. Jos\'e Gonz\'alez was supported by a grant from the Simons Foundation (Award Number 710443). Patricio Gallardo was partially supported by the National Science Foundation under Grant No. DMS-2316749. Javier Gonz\'alez Anaya is also grateful to Harvey Mudd College and Santa Clara University for their supportive environment, where part of this research was conducted.
We thank Alicia Lamarche for bringing to our attention the results that led to Corollary~\ref{cor:DC}. We are grateful to the referee for carefully reading the manuscript and providing valuable feedback, which has helped us improve the article overall.
\vspace{-2mm}

\tableofcontents


\section{Background}
\label{sec:Preliminarytdn}

In this section we provide the required background and notation for working on our moduli spaces.  Let $(z_1, \ldots, z_n)$ be a collection of $n$ distinct labeled points in $\mathbb{A}^d$. By considering the embedding
$\mathbb{A}^d \hookrightarrow \mathbb{P}^d$ given by $z_i \mapsto q_i:=[1:z_i]$, we obtain a configuration of $n$ distinct 
 labeled points in $\mathbb{P}^d$, all of them away from the hyperplane $H:=V(x_0)$.

We say two configurations of points $\{ q_1, \ldots, q_n \}$ and $\{ q_1', \ldots, q_n' \}$ are equivalent if there is an element $g\in\SL_{d+1}$ such that $g \cdot q_i = q_i'$ for all $i$ and $g|_{H}=id_{H}$. Let $G_d$ be the subgroup of $\SL_{d+1}$ that restricts to the identity in $H$. This is,
\begin{align*}
G_d:=
\{ 
g \in \SL_{d+1} \; | \; g|_{H} = id_H 
\}.
\end{align*}

Considering the diagonal $G_d$-action on $(\PP^d)^n$, these $G_d$-equivalence classes of point configurations are parameterized by the geometric quotient 

\begin{equation*}\label{eq:OpenTdn}
T_{d,n}^{\circ} = \mathcal{U}_{d,n} /  G_d, \qquad\text{where}\qquad \mathcal{U}_{d,n}
:=\{
(q_1, \ldots, q_n) \in (\mathbb{P}^d)^n
\; | \;
q_i \neq q_j, \; q_i \notin H, 
\text{ for all $i,j$}
\}.
\end{equation*}
From our choice of $H$, it follows that an element $g\in G_d$ is a matrix of the form
 \[
 g = 
  \begin{pmatrix}
   t^{-d} & 0 & \cdots & 0 \\
   s_{1} & t & \cdots & 0 \\
   \vdots & & \ddots & \\
   s_{d} & 0 & \cdots & t
  \end{pmatrix}\in\SL_{d+1},
 \]
and its action on a point in $\mathbb{A}^d = \mathbb{P}^d \setminus H$ is given by
\begin{align}\label{eq:actionAd}
g \cdot [1:z_1: \ldots : z_d] &= [t^{-d}:tz_1 + s_1: \ldots:  tz_n + s_d].   
\end{align} 
Therefore, $G_d$-orbits of points in $(\Aff^d)^n=(\mathbb{P}^d \setminus H)^n\subseteq(\mathbb{P}^d)^n$ are often referred to as a \emph{configuration of $n$ labeled points in $\Aff^d$ up to translation and scaling (or up to translation and homothety)}; see \cite{Chen-Gibney-Krashen,gallardo2018modular}.

In particular, $T_{1,n}^{\circ} \cong M_{0,n+1}$ because we can interpret the additional point in $\mathbb{P}^1$ as a hyperplane at infinity, and $G_1$ is the group of projective transformations that fix this point. Next, we describe a compactification that plays a central role in our work. 
\begin{lemma}[{\cite[Lemma 4.13]{GR}}]
\label{lemma:PPCompactification}
Fix a pair of integers $d,n\geq 2$, then we have an inclusion
\begin{align*}
  T^{\circ}_{d,n} \subseteq 
  \mathbb{P}^{d(n-1)-1}
  \cong 
  \left(  (\mathbb{A}^{n-1})^{d}  \setminus \vec{0} \right)
  / 
  \Gm.
\end{align*}
Moreover, we can choose a system of projective coordinates
so that the point 
\[
[x_{11}: x_{12}: \ldots : x_{1d} : x_{21}: x_{22} : \ldots x_{2d}:
 \ldots : x_{(n-1)1}: x_{(n-1)2} : \ldots : x_{(n-1)d}]\in \mathbb{P}^{d(n-1)-1}
\]
parameterizes the $G_d$-equivalence class associated to the collection of the following $n$ points in $\mathbb{P}^d$: 
\begin{align*}
q_1:=[1:x_{11}: \ldots : x_{1d}], & &
\ldots, & &
q_{n-1}:=[1:x_{(n-1)1}: x_{(n-1)2}: \ldots: x_{(n-1)d}], & &
q_n:= [1:0: \ldots: 0]
\end{align*}
\end{lemma}
The lemma above implies that the projective space $\mathbb{P}^{d(n-1)-1}$ can be identified with the moduli space of equivalence classes of $n$ points in $\mathbb{P}^d$ that lie away from a fixed hyperplane $H \subset \mathbb{P}^d$, up to the action of the subgroup $G_d \subseteq\SL_{d+1}$. In particular, the loci parameterizing configuration of colliding points can be described as follows. For each $I\subsetneq \{1,\dots,n\}$ with $2 \leq |I| \leq n-1$, we define subvarieties $\delta_{d,I}$ of $ \PP^{d(n-1)-1}$ parameterizing $G_d$-orbits  where the components indexed by $I$ coincide, this is, configurations in which the points $q_i\in\mathbb{P}^{d(n-1)-1}$ for $i \in I$ coincide. They are given by the equations 

 \begin{equation}\label{deltas}
 \delta_{d,I} =\begin{cases} 
\bigcap_{i,j \in I}V((x_{i1}-x_{j1},\ldots , x_{id}-x_{jd})) & \ \text{if } n\notin I;
\\
\bigcap_{i\in I\setminus n }V(( x_{i1} ,\ldots ,x_{id}))     & \  \text{if } n\in I.  
\end{cases}
\end{equation}

\begin{definition}\label{sec:admissibleWeight}
Consider an integer $n\geq 2$. 
We will refer to the set
\begin{align*}
\mathcal D_{n}^T :=
\left\{
(a_1, \ldots, a_n) \in \mathbb{Q}^{n} \; | \;
0 < a_i \leq 1 \textnormal{ for } i=1,\dots, n  \text{ and } 1 < a_1 + \ldots + a_n 
\right\} 
\end{align*}
as the \emph{domain of admissible weights}. Any $n$-tuple $\mathcal{A}\in \mathcal{D}_{n}^T $ will be referred to as a \emph{set of weights}. 
\end{definition}

\begin{remark} 
    The superscript $T$ in $ \mathcal{D}^T_n$ serves the purpose of distinguishing it from other domains of admissible weights that appear elsewhere in the literature in the context of other moduli problems.
\end{remark}

Let us define the set
\[
\wbuildtdn:=
\left\{
\delta_{d,I} \subset \PP^{d(n-1)-1} \ \left| \ \;
 \sum_{i \in I}a_i >1  \right.
\right\}.
\]

\begin{definition} \label{deftdno}
Let $\mathcal{A}\in \mathcal{D}^T_{n}$. 
A $G_d$-orbit in $(\mathbb{P}^d )^n$
parametrized by the open subvariety
\begin{align*}
\wtdno:=\bP^{d(n-1)-1} \setminus
\bigcup_{\delta_{d,I}\in \wbuildtdn} \delta_{d,I}.
\end{align*}
is called an \emph{$\mathcal{A}$-weighted configuration of $n$ labeled points in $\Aff^d$ up to translation and scaling}; 
 see~\eqref{eq:actionAd} and discussion afterward.
 For any $\mathcal{A}\in\mathcal{D}^T_n$ we have that $ T_{d,n}^\circ\subseteq\wtdno$, and this inclusion is often strict. 
\end{definition}

\begin{theorem}[{\cite[Theorem~1.4, Lemma~4.18]{GR}}] \label{theorem.wtdn}
For each $\mathcal{A} \in \mathcal D^T_{n}$, there exists a smooth projective variety $\wtdn$ of dimension $dn-d-1$
such that the following hold
\begin{enumerate}[label=(\roman*),ref=\thetheorem.(\roman*)]
\item\label{theorem.wtdn I} $\wtdn$ can be obtained as the blow-up of $\PP^{d(n-1)-1}$ along the subscheme $\bigcup_{\delta_{d,I}\in \wbuildtdn} \delta_{d,I}$,
 with scheme structure given by the product of the ideal sheaves of the $\delta_{d,I}$ in $\PP^{d(n-1)-1}$.

\item\label{theorem.wtdn II}
The boundary $\wtdn\setminus \wtdno$ is the union of $|\mathcal{H}_{\mathcal{A}}|$ 
smooth irreducible divisors. 

\item\label{theorem.wtdn III}
Any set of boundary divisors intersects transversely.
\end{enumerate}
\end{theorem}

The variety ${T_{d,n}^{\mathcal{A}}}$ is called the space of \emph{weighted pointed stable rooted trees} with respect to $\mathcal{A}$. 

\begin{remark}
    In the case where $\mathcal{A}=(1, \ldots, 1)$, the variety  $T^{\mathcal{A}}_{d,n}$ coincides with the Chen-Gibney-Krashen compactification  $T_{d,n}$ \cite{Chen-Gibney-Krashen} and we simply denote it by $T_{d,n}$ for ease of notation (see \cite[Definition 4.10]{GR}).
\end{remark}

\subsection{The toric compactification of the moduli space of $n$ points in affine space}\label{sec:PreliminaryToricModel}

Next we describe the compactification $T_{d,n}^{LM}$ of the space $T^{\circ}_{d,n}$ that is the focus of our work in this article. 
 This space generalizes the toric compactification of $M_{0,n}$ for $n \geq 3$ commonly known as the Losev-Manin space \cite{losev2000new-moduli-spaces}. 
Let ${T_{d,n}^{LM}}:=\wtdn$ be the space of weighted pointed stable rooted trees associated with any choice of weights $\mathcal{A}=(\epsilon_1,\dots,\epsilon_{n-1},1)\in \mathcal{D}^T_{n}$ such that $\sum\epsilon_i \leq 1$ (notice that $\epsilon_i > 0$ by the definition of $\mathcal{D}^T_{n}$).

For any such $\mathcal{A}$ we get the same locus $\mathcal{H}_{\mathcal{A}}$, and hence the same compactification $\wtdn$.
Moreover, for any such $\mathcal{A}$ the locus $\mathcal{H}_{\mathcal{A}}$ consists of subvarieties that are invariant under the standard toric variety structure of $\mathbb{P}^{d(n-1)-1}$ induced by its coordinate system from Lemma~\ref{lemma:PPCompactification}. 
Therefore, the compactification $T_{d,n}^{LM}$ is a toric variety. 

For each $1 \leq i \leq n-1$ and $1 \leq k \leq d$, 
let $e^k_i$ be the vector in $\mathbb{Z}^{d(n-1)} \subseteq \mathbb{R}^{d(n-1)}$ 
whose unique nonzero entry is a $1$ in the $d(i-1)+k$ position 
and let 
$\overline{e}^k_i$ denote the image of $e^k_i$ in the quotient $\mathbb{Z}^{d(n-1)} \! \! \! \! \left/ \left( \sum_{i=1}^{n-1} \sum_{k=1}^{d} e^k_i =0 \right). \right.$
For example, for $d=2$ and $n=3$, we have 
\begin{align*}
e^{1}_1={(1,0,0,0)}, \ e^{2}_1={(0,1,0,0)}, \ 
e^{1}_2={(0,0,1,0)},  \ e^{2}_2={(0,0,0,1)} \
\in \ 
\mathbb{Z}^4. 
\end{align*}
The $\overline{e}^k_i$ generate the rays of a fan in $\mathbb{R}^{d(n-1)} \! \! \! \! \left/ \left( \sum_{i=1}^{n-1} \sum_{k=1}^{d} e^k_i =0 \right) \right.$ describing $\mathbb{P}^{d(n-1)-1}$ as a toric variety. From the construction of $T_{d,n}^{LM}$ in Theorem~\ref{theorem.wtdn} as a blow-up of $\mathbb{P}^{d(n-1)-1}$ one obtains an explicit description of the rays of the fan of $T_{d,n}^{LM}$ as follows. 
\begin{corollary}[{\cite[Corollary 5.6]{GR}}]\label{toricTdn}
$T_{d,n}^{LM}$ is a toric variety and one can choose its fan to have rays generated by the vectors 
\begin{eqnarray*}
& \left\{ \left. \overline{e}^{k}_{i} \in \mathbb{Z}^{d(n-1)} \! / ( \sum_{i=1}^{n-1} \sum_{k=1}^{d} e^k_i =0 ) \, \right| \, 1 \leq i \leq n-1, \ 1 \leq k \leq d \right\}
\textnormal{\quad and }
\\
& \left\{
\left.
\sum_{i \in I} \left( \overline{e}^1_i+ \ldots+ \overline{e}^d_i \right) 
\in \mathbb{Z}^{d(n-1)} \! \! \! \! \left/ \! \! \left( \sum_{i=1}^{n-1} \sum_{k=1}^{d} e^k_i =0 \right) \right.
\ \right| \
1 \leq |I| \leq n-2, \ 
I \subsetneq \{1, \ldots, n-1 \}
\right\}.
\end{eqnarray*}
 In particular, $T_{1,n}^{LM}$ coincides with the Losev-Manin space $\overline{M}_{0,n+1}^{LM}$ introduced in \cite{losev2000new-moduli-spaces}.
\end{corollary}


\section{Higher-dimensional Losev-Manin spaces as toric Chow quotients}
\label{sec:AsToricChowQuotient}

In this section we prove Theorem~\ref{Main thm: I}. The Chow variety $\operatorname{Chow}((\PP^{d})^{n},\gamma)$ is a projective variety parameterizing formal sums with nonnegative integral coefficients of algebraic subvarieties of $(\PP^d)^{n}$ having homology class $\gamma$. 
The notion was introduced in \cite{kapranov1991quotients} for toric varieties, and in general in \cite{Kapranov-chow}. A detailed account of the general theory can be found in \cite[Chapter~1]{Kollar-chow}. 
In general, the Chow variety has a universal family parametrizing well-defined families of nonnegative, proper, algebraic cycles of fixed dimension and degree \cite[Theorem~3.21]{Kollar-chow}, and it can have arbitrarily bad singularities \cite[Theorem~1.1~M4]{vakil2006murphy}. Furthermore, there is a morphism from the seminormalization of the Hilbert scheme to the Chow variety \cite[Theorem~6.3]{Kollar-chow}.
For projective space $\mathbb{P}^d$, we have 
$H_*\left(\mathbb{P}^d,\mathbb{Z}\right) \cong \mathbb{Z}[x] /\left\langle x^{d+1}\right\rangle$, where the homology class of a subvariety $V$ is given by $\deg(V) x^{\dim(V)}$, where $\deg(V)$ is the degree of $V$. 
In particular, the homology class of a subvariety is uniquely determined by its degree and dimension. The homology of $(\PP^{d})^{n}$ is then fully described via K\"unneth's theorem.

We divide our presentation into four parts. First, we prove that there is an injective map  $T_{d,n}^{\circ}\hookrightarrow\operatorname{Chow}((\PP^{d})^{n-1},\gamma)$ for an appropriately chosen class $\gamma$; our Chow quotient $(\PP^d)^{n-1}\quotient_{Ch}\Gm$ is defined as the closure of the image of $T_{d,n}^{\circ}$ under this map. In the next subsection, we discuss the modular interpretation of the boundary strata of $T_{d,n}^{LM}$. Then, in the third subsection we show that the rational map $T_{d,n}^{LM}\dashrightarrow\operatorname{Chow}((\PP^{d})^{n-1},\gamma)$ given by the aforementioned injective map extends to a regular map. In the last subsection we prove that this extension defines an isomorphism from $T_{d,n}^{LM}$ to the normalization of our Chow quotient. 

\subsection{Definition of the Chow quotient}\label{chow: subsection 1} In Section~\ref{sec:Preliminarytdn} we defined the space $T_{d,n}^{\circ}=\mathcal{U}_{d,n} / G_d$, where 
\[
    \mathcal{U}_{d,n}=\{(q_1, \ldots, q_n) \in (\mathbb{P}^d)^n\; | \;q_i \neq q_j \text{ and } q_i \notin H \text{ for all $i,j$}\}
\]
and $G_d = \{g\in\SL_{d+1}\;|\; g|_H=\operatorname{id}_H\}$ acts diagonally on $(\PP^{d})^{n}$. In words, $T_{d,n}^{\circ}$ is the open subset of $T_{d,n}^{LM}$ parameterizing classes of $n$ distinct labeled points in $\mathbb{P}^d$, all of them away from the hyperplane $H=V(x_0)$.

By our choice of $H$, an element $g\in G_d$ is a matrix of the form
 \[
  \begin{pmatrix}
   t^{-d} & 0 & \cdots & 0 \\
   s_{1} & t & \cdots & 0 \\
   \vdots & & \ddots & \\
   s_{d} & 0 & \cdots & t
  \end{pmatrix}\in\SL_{d+1}.
 \]

\begin{definition}\label{chow: def diagonal torus}
    From this point forward we refer to the subgroup 
    \[
    \left\{ \left. \operatorname{diag}(t^{-d},t,\dots,t) \  \right|  \ t \in \Gm  \right\} \subseteq G_d
    \] 
    as the \emph{diagonal torus of $G_d$}, since it is clearly isomorphic to $\Gm$.
\end{definition}

\begin{proposition}\label{chow prop: iso quotient}
  There exists an open subset $\mathcal{V}\subseteq\mathcal{U}_{d,n-1}$ such that $\psi:\mathcal{V} \to T_{d,n}^{\circ}$ is a geometric quotient for the action of the diagonal torus $\Gm$.
\end{proposition}

Informally, this result says that, by translating the $n$\textsuperscript{th} point to the origin, it is possible to think of points in $T_{d,n}^{\circ}$ as parameterizing collections of $n-1$ distinct points in $\mathbb{A}^d\setminus\vec{0}$ up to scaling.

\begin{proof}
   Define the rational map $\mathcal{U}_{d,n-1}\dashrightarrow T_{d,n}^{\circ}$ mapping $(p_{1},\dots,p_{n-1})$ to $[(p_{1},\dots,p_{n-1},e_0)]_{G_{d}}$, where $e_0=[1:0:\cdots:0]\in\PP^{d}$. Its indeterminacy locus consists of points such that $p_i=e_0$ for some $i=1,\dots,n-1$. Let $\mathcal{V}$ be its complement and $\psi$ the restriction of this map to $\mathcal{V}$. We claim that $\psi$ is surjective and that its fibers are orbits of the $\Gm$-action. Since $T_{d,n}^{\circ}$ is normal, it follows that $\psi:\mathcal{V} \to T_{d,n}^{\circ}$ is a geometric quotient by~\cite[Remark 6.1.2]{dolgachev2003lectures}.
   
   First note that every class $[(q'_{1},\dots,q'_{n})]_{G_{d}}$ has a representative of the form $[(q_{1},\dots,q_{n-1},e_0)]_{G_{d}}$, because $q'_{n}\not\in H$. This implies that $\psi$ is surjective. 
  
  An element $(p_1,\dots,p_{n-1})\in\mathcal{V}$ maps to a class $[(q_{1},\dots,q_{n-1},e_0)]_{G_{d}}$ if and only if there exists some $g\in G_d$ such that $g\cdot p_i=q_i$ and $g\cdot e_0=e_0$. The latter condition implies that $g$ is an element of the diagonal torus $\Gm\subset G_d$. Hence, the fiber of $[(q_{1},\dots,q_{n-1},e_0)]_{G_{d}}$ is the $\Gm$-orbit of $(p_1,\dots,p_{n-1})$. 
\end{proof}

In fact, the fibers of $\psi$ are isomorphic to $\Gm$. Indeed, a point in $\PP^d$ is fixed by the diagonal torus if and only if it lies in $H$ or equals $e_0$. Thus, points in $\mathcal{V}$ have finite stabilizers and the claim follows. 

In order to define our Chow quotient, we have to show that the orbit closures $\overline{\Gm\cdot p}\subseteq(\PP^{d})^{n-1}$ have the same homology class $\gamma\in H_{2}((\PP^{d})^{n-1},\ZZ)$ for all $p\in\mathcal{V}$. Let us determine the cycle class of these orbit closures.

By K\"unneth's formula, there is a decomposition
\[
  H_{*}((\PP^{d})^{n-1},\ZZ) \cong \bigotimes_{i=1}^{n-1}H_{*}(\PP^{d},\ZZ),
\]
so every homology class is a linear combination of tensors composed of linear subspaces $[\PP^{m_1}]\otimes\cdots\otimes[\PP^{m_{n-1}}]$. Moreover, since every $\Gm$-orbit is one-dimensional, we are only interested in terms such that $\sum_{i=1}^{n-1} m_{i}=1$.

\begin{proposition}\label{prop: chow-class}
  Let $\mathcal{V}$ be as in Proposition~\ref{chow prop: iso quotient}. For any $p\in\mathcal{V}$,
  \[
    [\overline{\Gm\cdot p}] = \gamma := \sum_{i=1}^{n-1}\gamma_{i}\in H_{2}((\PP^{d})^{n-1},\ZZ),
  \]
  where $\gamma_{i}:=[\PP^{m_1}]\otimes\cdots\otimes[\PP^{m_{n-1}}]$ is the cycle class with $m_{i}=1$ and $m_{j}=0$ otherwise. 
\end{proposition} 
\begin{proof}
  We follow Kapranov's proof of the equivalent result; see \cite[Proposition~2.1.7]{Kapranov-chow}. 

  The $[\PP^{m_1}]\otimes\cdots\otimes[\PP^{m_{n-1}}]$ coefficient of $[\overline{\Gm\cdot p}]$ is the intersection number of the subvariety $\overline{\Gm\cdot p}$ with the product of general linear subspaces $L_{i}\subseteq\PP^{d}$ of codimension $m_{i}$. In the case at hand, $m_i=1$ while $m_j=0$ for all $j\neq i$. 
  
  Without loss of generality we compute the $\gamma_1$ coefficient, all other cases being identical. This coefficient is the intersection number between $\overline{\Gm\cdot p}$ and a generic hypersurface $L\times\PP^d\times\cdots\PP^d$, where $L\subseteq\PP^d$ is a hyperplane. The unique condition determining this intersection is the one imposed by $L$, so the coefficient of $\gamma_{1}$ can be computed by determining the cardinality of $\overline{\Gm\cdot p_1}\cap L$, where $p=(p_1,\dots,p_{n-1})\in(\PP^d)^{n-1}$.
  Let $p_1=[a_0:\cdots:a_d]$, then it can be easily verified that $\overline{\Gm\cdot p_1}$ is the line generated by $e_0$ and $[0:a_1:\dots:a_d]$. Note that these two points are precisely the boundary of the open orbit:
    \[
        \lim_{t\to 0} t\cdot p_1 
                                 = \lim_{t\to 0} [a_0:t^{d+1}a_1:\cdots:t^{d+1}a_d]
                                 = e_0,
    \]
    and
    \[
        \lim_{t\to\infty} t\cdot p_1 
                                 = \lim_{t\to\infty} [t^{-d-1}a_0:a_1:\cdots:a_d]
                                 = [0:a_1:\cdots:a_d].
    \]
    Since the orbit closure and $L$ are linear of complementary dimension, their intersection will consist of a single point as long as $L$ does not contain the orbit closure. If $L=V\left(\sum_{i=0}^d c_ix_i\right)$, then $\overline{\Gm\cdot p_1}\subseteq L$ if and only if $c_0=0$ and $\sum_{i=1}^dc_ia_i= 0$. In other words, a general hyperplane $L\in(\PP^d)^\vee \setminus V(c_0)\cap V\left(\sum_{i=1}^d c_ia_i\right)$ intersects the orbit closure at exactly one point.    
\end{proof}

The previous proposition allows us to define a set-theoretic map
\[
    \rho^{\circ}:T_{d,n}^{\circ}\hookrightarrow\operatorname{Chow}((\PP^{d})^{n-1},\gamma),
\]
mapping $x$ to the cycle $\overline{\psi^{-1}(x)}$ in the Chow variety. This function is clearly injective, owing to the fact that any two different points in $\mathcal{V}$ generate different $\Gm$-orbits in $(\PP^d)^{n-1}$.

In fact, the function $\rho^{\circ}:T_{d,n}^{\circ} \rightarrow \operatorname{Chow}((\mathbb{P}^{d})^{n-1},\gamma)$ is a morphism. 
Indeed, using the previous notation, let $\mathcal{X} \subseteq (\mathbb{P}^{d})^{n-1} \times T_{d,n}^{\circ}$ be the closure of the set 
\[
\left\{ \left. (p,x) \in (\mathbb{P}^{d})^{n-1} \times T_{d,n}^{\circ} \, \right| \, p \in \mathcal{V}  \textnormal{ and }  \psi(p)=x \right\},
\]
and let $\pi_2:\mathcal{X} \rightarrow T_{d,n}^{\circ}$ be the second projection restricted to $\mathcal{X}$. 
By construction, for each $x \in T_{d,n}^{\circ}$, we have that $\rho^{\circ}(x)$ is the class of the variety $\pi_2^{-1}(x)$ in $\operatorname{Chow}((\mathbb{P}^{d})^{n-1},\gamma)$. The morphism $\pi_2$ is proper, pure-dimensional and it has generically reduced fibers. 
Then, $\pi_2:\mathcal{X} \rightarrow T_{d,n}^{\circ}$ is a Chow-Cayley family in the sense of \cite[Section I.3]{Kollar-chow}. 
The Chow variety is a fine moduli space with a universal family for Chow-Cayley families over a seminormal base \cite[Section I.3]{Kollar-chow}.   
Therefore, $\rho^{\circ}:T_{d,n}^{\circ} \rightarrow \operatorname{Chow}((\mathbb{P}^{d})^{n-1},\gamma)$ is a morphism.

\begin{definition}\label{chow: definition chow quotient}
The Chow quotient $(\PP^{d})^{n-1}\quotient_{Ch}\Gm$ is the closure of the image of $T_{d,n}^{\circ}$ under $\rho^{\circ}$.
\end{definition}

This Chow quotient is a projective variety because it is an irreducible, closed subset of the Chow variety, which is a projective variety itself.

\subsection{The boundary strata of $T_{d,n}^{LM}$}\label{subsection: chow boundary} 
We remind the reader that a rooted tree graph is a connected acyclic graph with a distinguished vertex, called the \emph{root} of the tree.  
Given a rooted tree graph $V$, one defines a poset structure $P_V$ on $V$ where the elements are the vertices of the graph, and for any two vertices $v \neq w$, one has $v < w$ if and only if  
the unique path from the root to $w$ passes through $v$.  
In particular, the root is the unique minimal element of $P_V$.  
We denote rooted tree graphs by $V$; if vertices $v,w\in V$ satisfy $v < w$, then we say that $w$ is a \emph{descendant} of $v$ and $v$ is an \emph{ancestor} of $w$; if the two vertices are adjacent, we say that $w$ is a \emph{daughter} of $v$ and $v$ is the \emph{parent} of $w$.

For any $\mathcal{A}\in\mathcal{D}_n^T$ (see Definition~\ref{sec:admissibleWeight}), the closed points of $T_{d,n}^{\mathcal{A}}$ parameterize isomorphism classes of so-called $\mathcal{A}$-weighted $n$-pointed stable rooted trees.
These were first introduced in \cite[Section~2]{Chen-Gibney-Krashen} in the case that all weights equal $1$. The case for general weights was studied in \cite[Section~2.2]{GR}.  
By the mentioned references,  $\mathcal{A}$-weighted $n$-pointed stable rooted trees satisfy the following properties: they are $n$-pointed, reduced, not necessarily integral, equidimensional schemes of finite type over $\CC$ with simple normal crossing singularities only.   The combinatorial structure of an $\mathcal{A}$-weighted $n$-pointed stable rooted tree is described by its \emph{dual graph},    
which in this case is 
the graph whose vertices correspond to the irreducible components, and with edges joining pairs of vertices whose corresponding components have a nonempty intersection; for the general definition of the dual graph of a simple normal crossing divisor see \cite[Definition~7]{kollar2014simple}. 
Due to their construction, the dual graphs of $\mathcal{A}$-weighted $n$-pointed stable rooted trees are always rooted tree graphs.

The choice of weights has direct implications on the geometry of the $n$-pointed stable rooted trees parameterized by $T_{d,n}^{\mathcal{A}}$. From this point forward we focus on the case of $T_{d,n}^{LM}$, which corresponds to $\mathcal{A}=(\epsilon_1,\dots,\epsilon_{n-1},1)\in \mathcal{D}^T_{n}$ such that $\sum\epsilon_i \leq 1$  (notice that $\epsilon_i > 0$ by the definition of $\mathcal{D}^T_{n}$).    
In the rest of this subsection, we describe the key properties of the pointed stable, rooted trees parameterized by $T_{d,n}^{LM}$.

\begin{figure}
    \centering
    \includegraphics[width=10cm]{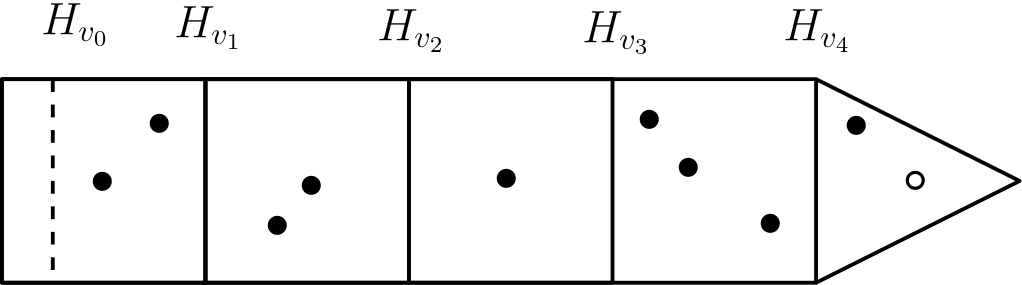}
\caption{A boundary point $\underline{X}\in T^{LM}_{2,10}$. Each square represents a component $X_v\cong\Bl_{e}\mathbb{P}^d$, with $v \neq v_{max}$, while the rightmost triangle represents the component $X_{v_{max}}\cong\mathbb{P}^d$. The common face between the $i$th and $(i+1)$th component represents the exceptional divisor $E_{v_i}$ in $X_{v_i}$ and the distinguished hyperplane $H_{v_{i+1}}$ in $X_{v_{i+1}}$. The black dots represent the labeled points with weights $<1$, while the hollow dot the point with weight $1$. 
}
\label{fig:deepest_stratum}
\end{figure}

\begin{lemma}[{\cite[Section~2.2]{GR}}]\label{chow lemma: rooted tree path}
Let $X$ be a weighted $n$-pointed stable rooted tree parameterized by $T_{d,n}^{LM}$. Then, $X$ satisfies the following properties; see Figure~\ref{fig:deepest_stratum}.
\begin{enumerate}
    \item The dual graph $V$ of $X$ is a rooted linear chain. The root vertex is denoted as $v_0$, while the maximal one as $v_{max}$. 

    \item Each irreducible component $X_v \subseteq X$ with $v < v_{max}$ is equipped with an isomorphism to the blow-up of $\PP^d$ at a point, we denote its exceptional divisor as $E_v$. The component $X_{v_{max}}$ is equipped with an isomorphism to $\PP^d$.  

    \item\label{chow: item hyperplane} Each $X_v\subseteq X$ is equipped with a fixed hyperplane $H_v\subseteq X_v$ disjoint from $E_v$.
    
    \item If $w$ is the daughter of $v$, then the intersection $X_v\cap X_w$ equals $E_v\subseteq X_v$ in the parent component, and $H_w \subseteq X_w$ in the daughter component.

    \item Each $X_v$ must have at least two different markings, all disjoint from $H_v$ and $E_v$; here a marking is either a marked point or $E_v$. Note that the hyperplane $H_v\subseteq X_v$ is not considered a marking of $X_v$.
    
    \item Any number of marked points in an irreducible component $X_v \subseteq X$ are allowed to collide provided that the sum of their weights remains $\leq 1$.

    \item The point of weight $1$ lies in $X_{v_{max}}$.
\end{enumerate}
\end{lemma}

From this point forward, for each $v\in V$ we will denote by $\overline{X}_v\cong\PP^d$ the image of the blow-down morphism from the previous lemma. Similarly, as we already did in item~\ref{chow: item hyperplane} above, we will not distinguish between the hyperplane $H_v\subseteq \overline{X}_v$ and its strict transform in $X_v$. 

\begin{definition}
    An $n$-pointed stable rooted tree is denoted as 
    \[
        \underline{X} = (\PP^{d-1}\hookrightarrow X = \cup_{v\in V}X_v;p_1,\ldots,p_n) \in T_{d,n}^{LM}.
    \] 
    Here $p_1,\dots,p_n$ are the marked points of $X$, the 
    dual graph $V$ of $X$ is a rooted tree graph, and $\PP^{d-1}  \hookrightarrow X$ is the hyperplane of the root component, called the \emph{root hyperplane} and denoted $H_{v_0}$. 
\end{definition}

This representation of a closed point $\underline{X}\in T_{d,n}^{LM}$ is defined only up to isomorphism. It will be convenient to work with specific representatives of such isomorphism classes, our choice is given in Lemma~\ref{Lemma rooted tree form}.

\begin{definition}[{\cite[Definition~2.0.4]{Chen-Gibney-Krashen}}]
\label{def:isomorphism_pointed_trees}
    Two $n$-pointed stable rooted trees $\underline{X}=(\PP^{d-1}\hookrightarrow X,p_1,\dots,p_n)$ and $\underline{Y}=(\PP^{d-1}\hookrightarrow Y,q_1,\dots,q_n)$ are isomorphic if there exists an isomorphism $f:X\to Y$ such that $f(p_i)=q_i$ for all $i=1,\dots, n$, and the following diagram commutes:
    \[
    \begin{tikzcd}[column sep=1.5em]
     & \PP^{d-1} \arrow[hookrightarrow]{dr}\arrow[hookrightarrow]{dl} \\
    X \arrow{rr}{f} && Y.
    \end{tikzcd}
    \]
\end{definition}

\begin{lemma} \label{Lemma rooted tree form}
    Every $n$-pointed stable rooted tree parameterized by $T_{d,n}^{LM}$ is isomorphic to one of the form $\underline{X}=(\PP^{d-1}\hookrightarrow X,p_1,\dots,p_{n-1},p_n)$ such that:
    \begin{enumerate}
        \item The last marked point is $p_n=e_0$ and, for all $v<v_{max}$, $X_v\cong\Bl_{e_0}\PP^d$;
        \item Every hyperplane $H_v\subseteq X_v$ is the strict transform of $V(x_0)\subseteq\overline{X}_v$.
        \item $p_i \neq p_n$ for all $1 \leq i \leq n-1$, while $p_i =p_j$ is allowed for any $1 \leq i,j \leq n-1  $. 
    \end{enumerate}
    Moreover, representatives of $\underline{X}$ satisfying properties (1) and (2) are uniquely defined up to the action of the diagonal torus $\Gm\subseteq G_d$ on each one of its irreducible components.
\end{lemma}
\begin{proof}
    Let $H \subseteq\mathbb{P}^d$ be a hyperplane and $p$ be a point disjoint from it. Then, there always exists an element $g \in \SL_{d+1}$ such that $g\cdot H = V(x_0)$ and $g\cdot p = e_0$. Indeed, the diagonal action of $\SL_{d+1}$ on $(\PP^d)^{d+2}$ is transitive on configurations of points in general position. This implies (1) and (2). 

    Next, let us prove the last claim. Consider two representatives $\underline{X}$ and $\underline{Y}$ of a closed point satisfying (1) and (2). 
    An isomorphism between $\underline{X}$ and $\underline{Y}$ induces a bijection between their components, mapping the root component to the root component. Moreover, it is easy to see inductively, starting from the root component, that all hyperplanes where components intersect are fixed pointwise.      
    In particular, the isomorphism restricts to the action of an element of $\Gm$ on the blow-down of each component, so the claim follows.
\end{proof}

\begin{conv}\label{chow: convention}
    From this point forward we consider all the objects parameterized by the closed points of $T_{d,n}^{LM}$ to be of the form presented in the previous lemma.
\end{conv}

\subsection{Extending the map}
In this section we show that $\rho^{\circ}$ extends uniquely to a regular map $\rho:T_{d,n}^{LM}\to\operatorname{Chow}((\PP^{d})^{n-1},\gamma)$. This is done by associating a cycle to each point configuration parameterized by $T_{d,n}^{LM}$. Note that there can be at most one such extension, since the image is dense in its closure and the Chow variety is separated. Moreover, the image of any such extension is contained in $(\PP^d)^{n-1}\quotient_{Ch}\Gm$, since the Chow quotient is closed in the Chow variety. 

Fix a closed point $\underline{X}\in T_{d,n}^{LM}$. By Lemma~\ref{chow lemma: rooted tree path}, each component $X_v\subseteq X$ is either $\PP^d$ or $\Bl_{e_0}\PP^d$. For any vertex $v\in V$, define the map $\varphi_v:X\to\overline{X}_v\cong\PP^d$ as follows. Consider any $p\in X_w$ such that $p=[a_0:\dots:a_d]\in\overline{X}_w$ and, if $w\neq v_{max}$, further assume that $p\neq e_0$. Then,
\begin{align}\label{chow: formula phi v}
    \varphi_{v}(p) = \left\{
    \begin{aligned}
        &[0:a_1:\dots:a_d],&\text{if } w<v;\\
        &[a_0:a_1:\dots:a_d],&\text{if }w=v;\\ 
        &e_0,&\text{if } v<w.
    \end{aligned}
    \right.
\end{align}

In fact, the maps $\varphi_v$ arise from the complete linear system of the line bundle $\cL_v := \varphi_v^*\cO_{\PP^d}(1)$. Indeed, a direct verification shows that
\begin{align}\label{chow: line bundle}
\cL_v|_{X_w} = \left\{
\begin{aligned}
    &\cO_{X_w}(H_w-E_w),&\text{if }w<v;\\
    &\cO_{X_w}(H_w),&\text{if }w=v;\\
    &\cO_{X_w},&\text{if }v<w.
\end{aligned}
\right.
\end{align}
Intuitively, points in an ancestor component of $X_v$ get inductively projected onto $H_v$, while points in its descendants get all contracted to $e_0$. See~\cite[Proposition~4.5]{gallardo2018modular} for a detailed account of the general properties of the line bundles $\cL_v$.

Two important observations are in order: First, note that the restriction $\varphi_v|_{X_w}$ is $\Gm$-invariant for all $w\neq v$, while it is $\Gm$-equivariant if $w=v$. On the other hand, note that the marked point $e_0\in X_{v_{max}}$ is mapped to $\varphi_v(e_0) = e_0$ for all $v\in V$. 

With this we are prepared to assign a homology class to any closed point in $T^{LM}_{d,n}$. 

\begin{definition}[{\cite[cf. Definition~4.7]{gallardo2018modular}}]
    Fix a closed point 
    \[
        \underline{X}=(\PP^{d-1} \hookrightarrow X = \cup_{v\in V}X_v,p_1,\ldots,p_{n-1},e_0) \in T_{d,n}^{LM}.
    \]
    For each $v\in V$, consider the point configuration $(\varphi_v(p_1),\dots,\varphi_v(p_{n-1}),e_{0})\in (\PP^d)^n$. The $v$-\emph{component configuration} is the point configuration
    \[
        \varphi_v(\underline{X}) = (\varphi_v(p_1),\dots, \varphi_v(p_{n-1}))\in(\PP^{d})^{n-1}.
    \]
    The \emph{configuration cycle} $Z(\underline{X})$ is the union of the $\Gm$-orbit closures of all component configurations:
    \[
        Z(\underline{X})= \bigcup_{v\in V}\overline{\Gm\cdot\varphi_v(\underline{X})}.
    \]   
\end{definition}
Let us show that the homology class of $Z(\underline{X})$ is the same as that one in Proposition~\ref{prop: chow-class}, independently of $\underline{X}$.

\begin{lemma}\label{chow: lemma J_v class}
    Let $\underline{X}\in T^{LM}_{d,n}$ be a closed point. For all $v\in V$, let $J_{v}=\{i : p_i\in X_{v}\}$ be the set of indices of the marked points in the component $X_{v}\subseteq X$. Then, 
    \[
        [\overline{\Gm\cdot\varphi_{v}(\underline{X})}] = \sum_{i\in J_v} \gamma_i,
    \]
    where $\gamma_{i}=[\PP^{m_1}]\otimes\cdots\otimes[\PP^{m_{n-1}}]$ is the cycle class with $m_{i}=1$ and $m_{j}=0$ otherwise.
\end{lemma}

\begin{proof}
    If $i\not\in J_v$, then $\varphi_{v}(p_i)\in\PP^d$ is a fixed point of the action, because any such marked point is either contained in $H$ or equal to $e_0$. These cases correspond to $p_i\in X_{w}$ for $w<v$ and $v<w$, respectively. Otherwise, if $i\in J_v$, then $\varphi_{v}(p_i)\in\PP^d$ generates a one-dimensional orbit $\Gm\cdot\varphi_{v}(p_i)\subseteq\PP^d$. The same argument as in the proof of Proposition~\ref{prop: chow-class} shows that the class of the orbit closure of $\varphi_v(\underline{X})=(\varphi_v(p_1),\dots,\varphi_v(p_{n-1}))$ is the sum $\sum \gamma_i$ over all indices $i$ such that $\varphi_v(p_i)$ has a one-dimensional orbit, this is, such that $i\in J_v$. The result follows.
\end{proof}

This simple lemma is quite helpful. For example, it immediately implies that the number of irreducible components of $Z(\underline{X})$ for $\underline{X}\in T_{d,n}^{LM}$ equals the number of vertices in its corresponding rooted tree graph. Moreover, it guarantees the extension of $\rho^{\circ}$ is viable:

\begin{proposition}\label{chow: proposition class boundary}
    The configuration cycle has the same homology class $[Z(\underline{X})] = \gamma$ for every closed point $\underline{X}\in T_{d,n}^{LM}$.
\end{proposition}
\begin{proof}
    Let $J_v$ be as in Lemma~\ref{chow: lemma J_v class}. Clearly $\{J_v\}_{v\in V}$ is a partition of $\{1,\dots,n-1\}$, so we have that
    \[
        [Z(\underline{X})]=\sum_{v\in V}[\overline{\Gm\cdot\varphi_{v}(\underline{X})}]=\sum_{v\in V} \sum_{i \in J_v} \gamma_i = \sum_{i=1}^{n-1}\gamma_i=\gamma.
    \]
\end{proof}

The following lemmas provide the key to extending $\rho^{\circ}$; see also \cite[Section 7.2]{Giansiracusa-Gillam}.

\begin{lemma}[{\cite[Lemma 3.18]{alexeev2023stable}}] \label{chow: alexeev}
    Let $X$ and $Y$ be proper varieties with $X$ normal. Let $\varphi: X \dashrightarrow Y$ be a rational map that is regular on an open dense subset $U \subseteq X$. Let $(C,0)$ be a regular curve and $f: C \longrightarrow X$ a morphism whose image intersects $U$. Let $g: C \longrightarrow Y$ be the unique extension of $f\circ\varphi$, which exists due to the properness of $Y$.

    Assume that for all $f$ with the same value of $f(0)$ there are only finitely many possibilities for $g(0)$. Then, $\varphi$ can be extended uniquely to a regular morphism $X \longrightarrow Y$.
\end{lemma}

The following result allows us to understand when two points are parameterized by the same cycle, even if the cycle is not irreducible. In the following statement, ``generic" refers to points with trivial stabilizer. However, this can be replaced with finite-dimensional stabilizers; see \cite[Remarks 3.2 and 3.14]{hu2005topological}.

\begin{lemma}\cite[Theorem 3.13]{hu2005topological}\label{yi hu lemma}
    Let $x$ and $y$ be two points in $X$ such that $\dim(G \cdot x) = \dim(G \cdot y)$. Then, the points $x$ and $y$ belong to the same cycle $Z_q$ parameterized for some 
    $q \in X\quotient_{Ch} G$   if and only if there is a generic holomorphic map from the complex unit disk $\phi: \Delta \to X$ with $\phi(0)=x$ and a holomorphic map $g: \Delta^* \to G$ from the punctured disk $ \Delta^*$ to $G$,  such that 
    \[
    y = \lim_{t \to 0} g(t) \cdot \phi(t).
    \]
\end{lemma}

\begin{lemma}\label{lemma:StrongLemma}
    Consider a stable rooted tree $\underline{X}_0\in T_{d,n}^{LM} \setminus T_{d,n}^{\circ}$ with dual graph $V$. Let $h:\Delta \to (\mathbb{P}^d)^{n-1}$ be an holomorphic map such that the diagram
\[
\begin{tikzcd}
    \Delta \ar[r,"h"] &  (\mathbb{P}^d)^{n-1} & \\
    \Delta^* \ar[r,"h^*"] \ar[u,hook] & \mathcal{V} \ar[r,"\psi"] \ar[u,hook] & T_{d,n}^{\circ}\subset T_{d,n}^{LM}
\end{tikzcd}
\]
commutes, and
\[
    \lim_{t \to 0} \psi ( h^*(t)) = \underline{X}_0.
\]
Then, for all vertices $v \in V$, there exists a holomorphic map $g_v(t):\Delta^*\to \{  \operatorname{diag}(1, t, \cdots, t)\,|\,t \in \mathbb{C}^*\}$ 
such that 

\vspace{-6mm}

\begin{align*}
 \lim_{t \to 0} g_v(t) \cdot h(t) =  \varphi_v(\underline{X}_0).
\end{align*}
\end{lemma}

\begin{proof}
Consider the configuration of points $h^*(t)=(q_1(t),\dots,q_{n-1}(t))$ with $q_i(t) = [1:x_{i1}(t): \ldots: x_{id}(t) ]$ (recall that by Convention~\ref{chow: convention} we set the $n$\textsuperscript{th} marked point of the configuration to be $q_n = e_0$). 

Write the Taylor expansion for the point $q_i(t)$ around $e_0\in\PP^d$ as 
\[
    x_{is}(t) = a_{is}t^{n_{is}} + O(n_{is}+1),
\]
where $O(n)$ denotes the terms of order $n$ and higher.

If $n_{is}\geq 1$ for all $s=1,\dots,d$, then blowing up the points colliding with $e_0$ reduces the exponents $n_{is}$ by one. Then, for each $v\in V$ there exist a number $n_v\in\NN$ such that 
$
\lim_{t \to 0} q_i(t) \in X_v
$
if and only if  $\min_{s}n_{is} = n_v$; see \cite[Section 1]{Fulton-MacPherson}. This allows to classify the points $q_i(t)$ into three groups:
\begin{align*}
q_i(t) = 
\begin{cases}
A_i(t), & \text{ if }    
\lim_{t \to 0} q_i(t) \in X_w,
\text{ with } w > v,
\\
B_i(t), & \text{ if }    
\lim_{t \to 0} q_i(t) \in X_v,
\\
C_i(t), & \text{ if }    
\lim_{t \to 0} q_i(t) \in X_w, 
\text{ with } w < v.
\end{cases}
\end{align*}
The Taylor expansion for each one of these cases is:
\begin{align*}
A_i(t) & = 
[1:a_{i1}t^{\alpha_{i1}} + O(\alpha_{i1}+1): \cdots: a_{id}t^{\alpha_{id}}+O(\alpha_{id}+1) ],
& &
\text{with all }\alpha_{is} > n_v,
\\
B_i(t) & = 
[1:a_{i1}t^{\beta_{i1}}+O(\beta_{i1}+1): \cdots: a_{id}t^{\beta_{id}}+O(\beta_{id}+1) ],
& &
\text{where }\beta_{is} \geq n_v,\;  \exists \beta_{is} = n_v,
\\
C_i(t) & = 
[1:a_{i1}t^{\gamma_{i1}}+O(\gamma_{i1}+1): \cdots: a_{id}t^{\gamma_{id}}+O(\gamma_{id}+1) ],
& &
\text{where }
\exists \gamma_{is} <  n_v.
\end{align*}
Next, we need to compute certain limits associated with the above Taylor expansion. Our choice of taking the $n$th marked point to be $e_0=[1:0:\dots:0]$ makes the following group convenient for our purposes. Define
\[
    D = \{ \operatorname{diag}(1,t,\dots,t)\,:\, t\in\Gm\}\subseteq \operatorname{GL}_{d+1}
\]
together with its standard action on $\mathbb{P}^d$, and its induced diagonal action on $(\mathbb{P}^d)^{n-1}$. The orbits of $D$ in $(\mathbb{P}^d)^{n-1}$ coincide set-theoretically with those of the diagonal torus $\Gm\subseteq G_d\subseteq\SL_{d+1}$ (see Definition~\ref{chow: def diagonal torus}). Moreover, the orbit closures (that is, the cycles) under both actions are the same. In particular, it is possible to apply Lemma~\ref{yi hu lemma} to the action of $G=D$, and the resulting conclusion will also hold for the action of the diagonal torus of $G_d$.

Now consider $g_v(t): \Delta^* \to D$ to be the matrix
\[
    g_v(t):=\text{diag}(1, t^{-n_v}, \ldots, t^{-n_v}).
\]
Then, it holds that 
\begin{align*}
\lim_{t \to 0} g_v(t) \cdot A_i(t) &= [1:0: \cdots: 0], 
\\
\lim_{t \to 0} g_v(t) \cdot B_i(t) &= 
[1:x_{i1}(0): \cdots: x_{id}(0)],   
\;  \text{ with } x_{is}(0) \neq 0 
\;  \text { for some } s
\\
\lim_{t \to 0} g_v(t) \cdot C_i(t) &= 
[0:x_{i1}(0): \cdots: x_{id}(0)], 
\;  \text{ with } x_{is}(0) \neq 0 
\;  \text { for some } s.
\end{align*}
This configuration of points is precisely $\varphi_v(\underline{X}_0)$ by construction.  
\end{proof}

With this we are now ready to construct the extension map.

\begin{proposition}\label{chow: proposition extension}
    The function $\rho:T_{d,n}^{LM}\to (\PP^d)^{n-1}\quotient_{Ch}\Gm$ given by $\rho(\underline{X})=Z(\underline{X})$ is a regular map.
\end{proposition}

\begin{proof}
    The restriction of $\rho$ to the interior $T_{d,n}^{\circ}$ is precisely $\rho^{\circ}$ from above. Our goal is to prove that this map extends uniquely to the boundary and that this extension is precisely $\rho$.

    Consider a curve $(C,0)\to T_{d,n}^{LM}$, where $(C,0)$ is as in Lemma~\ref{chow: alexeev}. Without loss of generality we may assume that $C\setminus 0$ maps to $T_{d,n}^{\circ}$ and $0$ maps to the boundary. Then, the configuration cycle $Z(\underline{X}_{\,s})$ of any point $s\in C\setminus 0$ has homology class $\gamma$, by Lemma~\ref{prop: chow-class}. Denote its limit in the Chow variety as $\lim_{s\to 0}Z(\underline{X}_{\,s})$. 
    
    By Lemma~\ref{chow: alexeev}, the existence of the extension will follow once we prove that $\lim_{s\to 0} Z(\underline{X}_{\,s})$ is uniquely determined by $\underline{X}_0$, the image of $0$ in this family. In fact, we will show that
    \[
        \lim_{s\to 0} Z(\underline{X}_{\,s}) = Z(\underline{X}_{\,0}).
    \]
    First, note that this equality will follow once we establish that $Z(\underline{X}_{\,0})\subseteq\lim_{s\to 0} Z(\underline{X}_{\,s})$. Indeed, define $D=\overline{\left(\lim_{s\to 0} Z(\underline{X}_{\,s})\right) \setminus Z(\underline{X}_{\,0})}$. Then, 
    \begin{align*}
        \gamma=[Z(\underline{X}_{\,0})] = \left[\lim_{s\to 0} Z(\underline{X}_{\,s})\right] + [D] = \gamma + [D] &\iff [D]=0.
    \end{align*}
     By a result of Kapranov \cite[Theorem~0.3.1]{Kapranov-chow}, each irreducible component of $\lim_{s\to 0} Z(\underline{X}_{\,s})$ is the $\Gm$-orbit closure of a point in $(\mathbb{P}^d)^{n-1}$
     \footnote[2]{We point out that \cite[Theorem~0.3.1]{Kapranov-chow} assumes that the stabilizer of a general point is trivial, but the proof provided there only uses that they are zero-dimensional, so it indeed applies in our present case. Alternatively; see the discussion before \cite[Proposition 3.6]{kapranov1991quotients}.}.
     Any such orbit closure has a nontrivial homology class, since it is either a point or a one-dimensional orbit closure, so $[D]=0$ if and only if $D=\emptyset$.

    Let us proceed to show that $Z(\underline{X}_{\,0})\subseteq\lim_{s\to 0} Z(\underline{X}_{\,s})$. It is enough to show that each irreducible component of the former is contained in the latter, this is, that 
    \[
        \overline{\Gm\cdot
        \varphi_v(\underline{X}_{\,0})}\subseteq\lim_{s\to 0} Z(\underline{X}_{\,s})\subseteq(\PP^d)^{n-1}
    \]
    for all vertices $v$ in the tree graph of $\underline{X}_{\,0}$. Furthermore, it is enough to prove that 
    \[
        \varphi_v(\underline{X}_{\,0})=(\varphi_v(p_1),\dots,\varphi_v(p_{n-1}))\in\lim_{s\to 0} Z(\underline{X}_{\,s})\subseteq(\PP^d)^{n-1}.
    \]
    This containment follows from Lemma~\ref{yi hu lemma} together with Lemma~\ref{lemma:StrongLemma}. Indeed, let $k:\Delta^* \to T_{d,n}^{\circ}$ be an holomorphic map obtained by restricting the map $(C,0)\to T_{d,n}^{LM}$ to an analytic neighborhood $\Delta$ of $0\in C$, so that $\lim_{t \to 0} k(t) = \underline{X}_{\,0}$.  
    Consider a holomorphic map $h^*:\Delta^* \to (\mathbb{P}^d)^{n-1}$ such that $k = \psi \circ h^*$, where $\psi:\mathcal{V}\to T_{d,n}^{\circ}$ is the geometric quotient from Proposition~\ref{chow prop: iso quotient}.  
    To see that such $h^*$ exists, it is enough to construct a section of $\psi$. Now, Lemma~\ref{lemma:PPCompactification} gives us an open embedding $T_{d,n}^{\circ} \subseteq \mathbb{P}^{d(n-1)-1}$ and a section of the rational map $(\mathbb{P}^d)^{n-1} \dashrightarrow \mathbb{P}^{d(n-1)-1}$, which together yield the desired section of $\psi$.

    By Lemma~\ref{lemma:StrongLemma}, there exists a $g_v:\Delta^* \to \{ t \in \mathbb{C}^* \, | \, \operatorname{diag}(1, t, \cdots, t)\}$ such that 
    \[
        \lim_{t \to 0} g_v(t) \cdot h^*(t) = \varphi_v(\underline{X}_0)\in (\mathbb{P}^d)^{n-1}.
    \]
    Notice that the action of $G:= \{ t \in \mathbb{C}^* \, | \, \operatorname{diag}(1, t, \cdots, t)\}$ on $(\mathbb{P}^d)^{n-1}$ and that of our diagonal torus $\mathbb{C}^*= \{ t \in \mathbb{C}^* \, | \, \operatorname{diag}(t^{-d}, t, \cdots, t)\}$ produce the same orbits, the same cycles and the same Chow quotient 
    $(\PP^d)^{n-1}\quotient_{Ch}G = (\PP^d)^{n-1}\quotient_{Ch}\Gm$. 
    By Lemma~\ref{yi hu lemma} 
    applied to the group $G$ and the Chow quotient 
    $(\PP^d)^{n-1}\quotient_{Ch}G = (\PP^d)^{n-1}\quotient_{Ch}\Gm$
    both $\lim_{t \to 0} g_v(t) \cdot h^*(t)$ and $\varphi_v(\underline{X}_0)$ belong to the same cycle for all $v$ in $V$. 
    
    On the other hand, $\lim_{t \to 0} g_v(t) \cdot h^*(t)\in\lim_{s \to 0} Z(\underline{X}_{\,s})$ because $g_v(t) \cdot h^*(t)\in Z(\underline{X}_{\,t})$ for all $t\in\Delta^*$. 
    To see this, consider the Chow-Cayley family $\mathcal{C}$ obtained from the map $C\to(\PP^d)^{n-1}\quotient_{Ch}\Gm$. Then, $\mathcal{C}\subseteq (\PP^d)^{n-1}\times C$ is such that $\pi_2^{-1}(s)$ is the cycle parameterized by the image of $s\in C$ in the Chow quotient. It follows that the limit $\lim_{t \to 0} g_v(t)\cdot h^*(t)\in\mathcal{C}$ lies in the cycle $\pi_2^{-1}(0)$, which is precisely $\lim_{s \to 0} Z(\underline{X}_{\,s})$.
    This implies that $\varphi_v(\underline{X}_0)\in\lim_{s \to 0} Z(\underline{X}_{\,s})$ for all $v\in V$.
\end{proof}

\subsection{Isomorphism with the normalization of the Chow quotient} To conclude this section we prove Theorem~\ref{Main thm: I}. Explicitly, we prove that the map $\rho:T_{d,n}^{LM}\to(\PP^d)^{n-1}\quotient_{Ch}\Gm$ from Proposition~\ref{chow: proposition extension} is an isomorphism. 

\begin{lemma}\label{chow: corollary injectivity}
    The map $\rho:T_{d,n}^{LM}\to(\PP^d)^{n-1}\quotient_{Ch}\Gm$ from Proposition~\ref{chow: proposition extension} is bijective. 
\end{lemma}

\begin{figure}
    \centering
    \includegraphics[width=\textwidth]{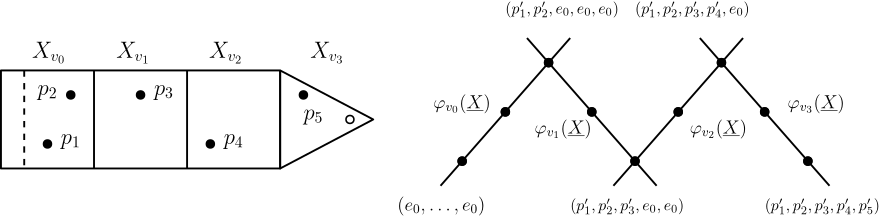}
\caption{The left-hand side figure represents a boundary point $\underline{X}\in T_{d,6}$, as in Figure~\ref{fig:deepest_stratum}. The right-hand side figure is its corresponding cycle $Z(\underline{X})\subseteq(\PP^d)^{5}$. The cycle $Z(\underline{X})$ has one irreducible component per vertex $v\in V$, given by $\overline{\Gm\cdot \varphi_{v}(\underline{X})}\cong\mathbb{P}^1$.
 Denote by $p'\in\mathbb{P}^d$ the projection of $p\in\PP^d$ from $e_0$ to the hyperplane $H=V(x_0)$. Each line in the right-hand side figure corresponds to an irreducible component of $Z(\underline{X})$. Each point $\varphi_{v_i}(\underline{X})\in(\mathbb{P}^d)^{5}$ is represented as a black dot. The other two black dots in each line segment represent the two limit points of the orbit $\Gm\cdot \varphi_{v_i}(\underline{X})$. } 
\label{chow: figure z(x)}
\end{figure}

\begin{proof}
Consider two closed points $\underline{X}$ and $\underline{Y}$ of $T_{d,n}^{LM}$. To prove injectivity, we show that if 
$Z(\underline{X})=Z(\underline{Y})$, then 
$\underline{X}$ and $\underline{Y}$ are isomorphic (Definition \ref{def:isomorphism_pointed_trees}). By Lemma~\ref{Lemma rooted tree form}, the points can be taken to be of the form $\underline{X}=(X,p_1,\dots,p_{n-1},e_0)$ and $\underline{Y}=(Y,q_1,\dots,q_{n-1},e_0)$, and showing they are isomorphic amounts to proving two conditions:
\begin{enumerate}
    \item 
    The underlying varieties $X$ of $\underline{X}$ and $Y$ of $\underline{Y}$ are isomorphic.
    By item 2 in Lemma~\ref{chow lemma: rooted tree path}, the dual graph of $X$, resp. $Y$, determines the variety $X$, resp. $Y$, uniquely up to isomorphism, because their components are isomorphic to either 
    $\mathbb{P}^d$ or $\Bl_{pt.}\mathbb{P}^{d-1}$. Therefore, to show $X$ and $Y$ are isomorphic it suffices to prove that the dual graphs of $X$ and $Y$ are the same.

    \item If $V$ denotes the dual graph of $X$ and $Y$, for each $v\in V$ there is an equality
    \[
        J_v:=\{i\,:\, p_i\in X_v\} = \{i\,:\, q_i\in Y_v\},
    \]
    and there exists some $t_v\in\Gm$ such that $p_i = t_v\cdot q_i\text{ for all } i\in J_v$. In words, the labeled points in each irreducible component of $X$ and $Y$ are the same up to the action of the diagonal torus $\Gm\subseteq G_d$. 
\end{enumerate}

Let us prove the first claim. The dual graph of any closed point of $T_{d,n}^{LM}$ is a path graph by definition. Therefore, the first condition amounts to proving that $Z(\underline{X})$ determines the number of vertices in the dual graph $V$ of $X$, and the same for $Z(\underline{Y})$ and $Y$. For each vertex $v\in V$, define $Z_v := \overline{\Gm\cdot\varphi_v(\underline{X})}$. By Lemma~\ref{chow: lemma J_v class}, we have $[Z_v] \neq [Z_w]$ for any two distinct vertices $v, w \in V$. In particular, this implies that $Z_v \neq Z_w$. 

Since the $Z_v$ are distinct irreducible components of $Z(\underline{X}) = \bigcup_{v \in V} Z_v$, we conclude that $Z(\underline{X})$ has the same number of irreducible components as $X$. The same reasoning applies to $\underline{Y}$. Therefore, if $Z(\underline{X}) = Z(\underline{Y})$, then $X$ and $Y$ have the same number of irreducible components arranged in a path configuration, so $X$ and $Y$ are isomorphic.

To show the second claim we begin by identifying the two limit points of the orbit $\Gm\cdot\varphi_v(\underline{X})$ for each $v\in V$. Given a point $p=[a_0:\dots:a_d]\in\mathbb{P}^d$, define $p'$ to be its projection to the root hyperplane $H=V(x_0)$, that is, $p' = [0:a_1:\dots:a_d]$. Then, 
\[
    \lim_{t\to 0}t\cdot\varphi_v(\underline{X}) = (r_1,\dots,r_{n-1})\quad\text{and}\quad \lim_{t\to\infty}t\cdot\varphi_v(\underline{X}) = (s_1,\dots,s_{n-1}),
\]
where, by Equation~\eqref{chow: formula phi v}, the individual components are
\begin{align*}
r_i = \lim_{t\to 0}t\cdot\varphi_v(p_i)=\left\{
    \begin{aligned}
    p_i', \text{ if } p_i\in X_w \text{ for } w<v;\\
    e_0, \text{ if } p_i\in X_w \text{ for } w=v;\\
    e_0, \text{ if } p_i\in X_w \text{ for } v<w;
    \end{aligned}
\right.
\quad\text{and}\quad
s_i = \lim_{t\to\infty}t\cdot\varphi_v(p_i)=\left\{
    \begin{aligned}
    p_i', \text{ if } p_i\in X_w \text{ for } w<v;\\
    p_i', \text{ if } p_i\in X_w \text{ for } w=v;\\
    e_0, \text{ if } p_i\in X_w \text{ for } v<w.
    \end{aligned}
\right.
\end{align*}

In particular, note that if $v<w$ and $w$ is the daughter of $v$, then $\lim_{t\to\infty} t\cdot\varphi_v(\underline{X})=\lim_{t\to 0}t\cdot\varphi_w(\underline{X})$. Moreover, in this case $\overline{\Gm\cdot\varphi_v(\underline{X})}$ and $\overline{\Gm\cdot\varphi_w(\underline{X})}$ are the only two irreducible components of $Z(\underline{X})$ passing through this point. See Figure~\ref{chow: figure z(x)}.

Let us now order the vertices of $V$ as $v_0<v_1<v_2<\cdots$, starting from the root vertex $v_0$. Then, by the previous formulas, we have that:
\begin{itemize}
    \item The orbit closure $Z_0 = \overline{\Gm\cdot\varphi_{v_0}(\underline{X})}$ is the unique irreducible component of $Z(\underline{X})$ containing the point $(e_0,\dots,e_0)$. Similarly, $Z_0'=\overline{\Gm\cdot\varphi_{v_0}(\underline{Y})}$ is the unique irreducible component of $Z(\underline{Y})$ containing the point $(e_0,\dots,e_0)$. Since $Z(\underline{X})=Z(\underline{Y})$, we conclude that $Z_0=Z_0'$, that is,
    \[
        \overline{\Gm\cdot\varphi_{v_0}(\underline{X})} = \overline{\Gm\cdot\varphi_{v_0}(\underline{Y})}.
    \]

    \item The orbit closure $Z_1 = \overline{\Gm\cdot\varphi_{v_1}(\underline{X})}$ is the unique irreducible component of $Z(\underline{X})$, other than $Z_0$, containing the point $\lim_{t\to\infty}t\cdot\varphi_{v_0}(\underline{X})$. Similarly, the orbit closure $Z_1' = \overline{\Gm\cdot\varphi_{v_1}(\underline{Y})}$ is the unique irreducible component of $Z(\underline{Y})$, other than $Z_0'$, containing the point $\lim_{t\to\infty}t\cdot\varphi_{v_0}(\underline{Y})$. However, since $Z_0=Z_0'$ and $Z_1, Z_1'\neq Z_0$, the equality $Z(\underline{X})=Z(\underline{Y})$ implies that $Z_1=Z_1'$, that is,
    \[
        \overline{\Gm\cdot\varphi_{v_1}(\underline{X})} = \overline{\Gm\cdot\varphi_{v_1}(\underline{Y})}.
    \]
    \end{itemize}
    \[ \vdots \]
    \begin{itemize}
    \item In general, for all $i=1,\dots,|V|$, the orbit closure $Z_i = \overline{\Gm\cdot\varphi_{v_i}(\underline{X})}$ is the unique irreducible component of $Z(\underline{X})$, other than $Z_{i-1}$, containing the point $\lim_{t\to\infty}t\cdot\varphi_{v_{i-1}}(\underline{X})$. Similarly, the orbit closure $Z_i' = \overline{\Gm\cdot\varphi_{v_i}(\underline{Y})}$ is the unique irreducible component of $Z(\underline{Y})$, other than $Z_{i-1}'$, containing the point $\lim_{t\to\infty}t\cdot\varphi_{v_{i-1}}(\underline{Y})$. Then, using the same argument as before, we obtain that
    \[
        \overline{\Gm\cdot\varphi_{v_i}(\underline{X})} = \overline{\Gm\cdot\varphi_{v_i}(\underline{Y})}.
    \]
\end{itemize}
The previous argument shows that $\overline{\Gm\cdot\varphi_{v}(\underline{X})} = \overline{\Gm\cdot\varphi_{v}(\underline{Y})}$ for all $v\in V$. For any fixed $v$, this implies that there exists some $t_v\in\Gm$ such that
\[
    (\varphi_{v}(p_1),\dots, \varphi_{v}(p_{n-1})) = 
    t_v\cdot (\varphi_{v}(q_1),\dots, \varphi_{v}(q_{n-1})).
\] 
Hence, we have that $p_i = t_v\cdot q_i$ for all $i\in J_v$, because $\varphi_v(p_i)=p_i$ and $\varphi_v(q_i)=q_i$ for all $i\in J_v$, by Equation~\eqref{chow: formula phi v}. This concludes the proof of injectivity.

Proving that $\rho$ is surjective is significantly simpler. The map is projective because both its domain and codomain are projective varieties. In addition, it is a dominant map by construction, so the claim follows.

\end{proof}

Let $\rho:T_{d,n}^{LM}\to(\PP^d)^{n-1}\quotient_{Ch}\Gm$ be the extension constructed in Proposition~\ref{chow: proposition extension} and $\pi:N\to (\PP^d)^{n-1}\quotient_{Ch}\Gm$ the normalization of the Chow quotient. Since $T_{d,n}^{LM}$ is smooth, the universal property of the normalization gives a factorization 
\[
    \rho=\pi\circ\tilde{\rho}:T_{d,n}^{LM} \stackrel{\tilde{\rho}}{\to} N \stackrel{\pi}{\to} (\PP^d)^{n-1}\quotient_{Ch}\Gm.
\]  
The main result of this section is the following.

\begin{theorem}[Theorem~\ref{Main thm: I}]\label{chow main thm}
    The map $\tilde{\rho}:T_{d,n}^{LM}\to N$ is an isomorphism.
\end{theorem}
\begin{proof}
    We use the following version of Zariski's Main Theorem \cite[Theorem~12.83]{GortzWedhorn}: If $f:X\to Y$ is a bijective proper morphism of varieties, with $Y$ normal, then $f$ is an isomorphism. 
    
    The result will then follow once we show that $\tilde{\rho}$ is a bijection. The map $\tilde{\rho}$ is injective and dominant, since $\rho$ is bijective by Lemma~\ref{chow: corollary injectivity}. The properness of $\tilde{\rho}$ then implies it is surjective.
    \end{proof}


\section{Mori dream space property of the Chen-Gibney-Krashen compactification \texorpdfstring{$T_{d,n}$}{Td,n}}\label{sec:NotMDS}

In this section we prove Theorem~\ref{Main thm: II}, Theorem~\ref{introduction: blowup theorem quotient} and Theorem~\ref{theorem.Tdn.not.MDS} from the introduction. 
First, in Section~\ref{subsection.Hausen.Suss} we review some results of Hausen and Su{\ss} \cite{hausen2010cox} about the Cox rings of varieties with torus actions. 
Then, in Section~\ref{subsection.blowup.subtorus.closure} we show Theorem~\ref{introduction: blowup theorem quotient} (Theorem~\ref{theorem.blow-up.torus.along.subtorus}) about the Cox ring of the blow-up of a toric variety along the closure of a subtorus. 
Finally, in Section~\ref{subsection.non-MDS.Tdn} we prove Theorem~\ref{Main thm: II} (Theorem~\ref{thm:MainLM.Part.b}) and Theorem~\ref{theorem.Tdn.not.MDS} (Theorem~\ref{theorem.Tdn.not.MDS.chapter.4}) regarding the Mori dream space property of the compactification $T_{d,n}$.

\subsection{The Cox ring of an algebraic variety with an effective torus action} \label{subsection.Hausen.Suss} Let $X$ be a normal algebraic variety with an effective action of a torus $T$.
Recall that an effective group action is one where the only element that acts trivially is the identity. 
The subset $X_0$ of $X$ consisting of the points with zero-dimensional stabilizer is nonempty, open and $T$-invariant.
By \cite{sumihiro1974equivariant}, there exists an irreducible normal but possibly non-separated prevariety $Y_0$ and a morphism $\pi: X_0 \rightarrow Y_0$ that is a geometric quotient for the induced action of $T$ on $X_0$. 
Following \cite[Section 1]{hausen2010cox}, we know that there exists a \emph{separation} of $Y_0$, that is, a variety $Y$ and a rational map $\theta: Y_0 \dashrightarrow Y$ 
which is defined on a big open subset $U \subseteq Y_0$ and maps $U$ locally isomorphically onto a big open subset $V \subseteq Y$.
Here, a big open subset means an open subset whose complement has codimension at least two.  
The geometric quotient $\pi: X_0 \rightarrow Y_0$ is unique, but the separation $\theta: Y_0 \dashrightarrow Y$ is not unique in general. Let us fix a separation $\theta: Y_0 \dashrightarrow Y$. 
In this setting, and further assuming that $X$ is complete and has a finitely generated class group,
\cite[Theorem 1.2]{hausen2010cox} says that the Cox ring of $X$ is finitely generated if and only if the Cox ring of $Y$ is finitely generated.
Moreover, the same theorem gives a presentation of the Cox ring of $X$ as an algebra over the Cox ring of $Y$ in terms of generators and relations.  

The next two examples will be useful in what follows. 

\begin{example}[Toric downgrades]  \label{example.toric.case}  
Let $X$ be a toric variety with torus $T_X$ associated to the fan $\Delta_X$ in the lattice $N_X$. 
Let $T$ be a subtorus of $T_X$ arising from an inclusion of lattices $N \hookrightarrow N_X$. 
Hence, $T \rightarrow T_X$ is a closed immersion and the quotient $N_X/N$ is torsion-free. 
Let us assume that the maximal cones in the fan $\Delta_X$ are one-dimensional. 
Let $\rho_1,\rho_2, \ldots, \rho_d$ be the rays in $\Delta_X$, with images in $N_X/N$ denoted by $\overline{\rho_1},\overline{\rho_2}, \ldots, \overline{\rho_d}$, respectively. 
We may assume the rays are ordered such that there exists $0 \leq r \leq d$ such that $\overline{\rho_i}$ is one-dimensional for $i \leq r$ and zero-dimensional for $i >r$. 
Let $\rho'_1,\rho'_2, \ldots, \rho'_s$ be the rays $\overline{\rho_1},\overline{\rho_2}, \ldots, \overline{\rho_r}$ listed in any order but without repetitions. 
As usual, we denote the affine toric variety associated to a cone $\sigma$ by $U_{\sigma}$. 
Let us review \cite[Remark 5.10]{hausen2010cox}, which explains what the geometric quotient $\pi: X_0 \rightarrow Y_0$ and a separation $\theta: Y_0 \dashrightarrow Y$ from the previous discussion are in the present setting. 

In the present setting, the open subset $X_0$ of $X$ consisting of the points with zero-dimensional stabilizer is precisely the toric subvariety corresponding to the subfan of $\Delta_X$ obtained by removing the rays $\rho_i$ for $i > r$.
The prevariety $Y_0$ is obtained by gluing the affine toric varieties $U_{\overline{\rho_1}},U_{\overline{\rho_2}}, \ldots, U_{\overline{\rho_r}}$ along their common open torus, and the geometric quotient $\pi: X_0 \rightarrow Y_0$
is obtained by gluing the maps $U_{\rho_i} \rightarrow U_{\overline{\rho_i}}$ induced by the projection $N_X \rightarrow N_X/N$. 
The variety $Y$ is the toric variety obtained by gluing the affine toric varieties $U_{{\rho_1'}},U_{{\rho_2'}}, \ldots, U_{{\rho_s'}}$ along their common open torus, and the separation $\theta: Y_0 \dashrightarrow Y$ is induced by gluing the identity maps $U_{\overline{\rho_i}} \rightarrow U_{\rho'_j}$ for all indices such that $\overline{\rho_i}=\rho'_j$. 
\end{example}

\begin{example}[Product with a variety with a trivial action]  \label{trivial.product}
The construction above commutes with taking a product with a factor that has a trivial action in the following sense. 
Suppose that $X$ is a normal variety with an effective action of a torus $T$. Let $X_0$, $Y_0$ and $Y$ be defined as above, with $\pi: X_0 \rightarrow Y_0$ denoting the geometric quotient and $\theta: Y_0 \dashrightarrow Y$ denoting a separation. 
Let $X':=X \times Z$ with the action of $T$ induced by the given action of $T$ on $X$ and a trivial action of $T$ on the factor $Z$. 
Then, the set of points $X'_0$ in $X'$ with zero-dimensional stabilizers is precisely $X'_0=X_0 \times Z$.
Moreover, $\pi \times \operatorname{Id}_{Z}: X_0 \times Z \rightarrow Y_0 \times Z$ is a geometric quotient for the $T$-action and $\theta \times \operatorname{Id}_{Z}: Y_0 \times Z \dashrightarrow Y \times Z$ is a separation. 
\end{example}

\subsection{The Cox ring of the blow-up of a toric variety along the closure of a subtorus of its torus}  \label{subsection.blowup.subtorus.closure}
In this section we prove Theorem~\ref{introduction: blowup theorem quotient}. 
Let $T_1$ be a torus and $T_2$ be a subtorus. 
We get an induced action of $T_2$ on $T_1$ and its subvariety $T_2$ is $T_2$-invariant. 
This further induces an action of $T_2$ on the blow-up $\operatorname{Bl}_{T_2}(T_1)$. 
The quotient homomorphism $T_1 \rightarrow T_1/T_2$ is a geometric quotient for the $T_2$-action and $T_1/T_2$ is a torus with identity element $T_2/T_2$. 
 
\begin{lemma} \label{lemma.blow-up.torus.along.subtorus}
  Let $T_1$ be a torus and $T_2$ be a subtorus. 
Then, there is a $T_2$-equivariant isomorphism 
\[
\operatorname{Bl}_{T_2}(T_1) \ \stackrel{\cong}{\longrightarrow} \ T_2 \times \operatorname{Bl}_{T_2/T_2}(T_1/T_2), 
\]
where the $T_2$-action on the right-hand side is trivial on the second factor. 
  In particular, the quotient $T_1 \rightarrow T_1/T_2$ induces a geometric quotient 
$\pi: \operatorname{Bl}_{T_2}(T_1) \rightarrow \operatorname{Bl}_{T_2/T_2}(T_1/T_2)$ for the induced $T_2$-action.
\end{lemma}

\begin{proof}
Let $N_1$ and $N_2$ be the lattices of one-parameter subgroups of $T_1$ and $T_2$. 
We can consider $N_2$ as a saturated sublattice of $N_1$. 
Let us choose a complement for $N_2$, that is, 
a saturated sublattice $N_3$ of $N_1$ such that $N_1 = N_2 \oplus N_3$.
This induces a decomposition $T_1=T_2 \times T_3$, where $T_3$ is the torus associated to the lattice $N_3$. 
We can identify $T_1/T_2$ with $T_3$ and the point $T_2/T_2$ with the identity element $t_3$ of $T_3$. 
The projection $f: T_2 \times T_3 \rightarrow T_3$ is a flat morphism and $f^{-1}(t_3)=T_2 \times \{t_3 \}$. Hence, by the commutativity of blow-ups with flat base change, we get a commutative diagram with cartesian squares
\[
\begin{tikzcd}
 \operatorname{Bl}_{T_2 \times \{t_3\}} (T_2 \times T_3) \arrow{d} \arrow{r} \dar\drar[phantom, "\square"]
  & T_2 \times T_3 \arrow{d} \arrow{r} \dar\drar[phantom, "\square"]
  & T_2  \arrow{d} \\
 \operatorname{Bl}_{t_3}T_3 \arrow{r} & T_3 \arrow{r}
& \operatorname{Spec}(k) \end{tikzcd}.
\]
Therefore, $\operatorname{Bl}_{T_2}(T_1)=\operatorname{Bl}_{T_2 \times \{t_3\}} (T_2 \times T_3)$ is isomorphic to $T_2 \times \operatorname{Bl}_{T_2/T_2}(T_1/T_2)=T_2 \times \operatorname{Bl}_{t_3}(T_3)$ and the isomorphism is $T_2$-equivariant if we let $T_2$ act trivially on $\operatorname{Bl}_{t_3}(T_3)$. 
This proves the first claim and from this the claim about the geometric quotient $\pi$ follows at once. 
\end{proof}

Now we prove Theorem~\ref{introduction: blowup theorem quotient}, which gives us a tool to study the finite generation of the Cox rings of toric varieties blown up along the closure of a subtorus of their torus.

\begin{theorem}[Theorem~\ref{introduction: blowup theorem quotient}] \label{theorem.blow-up.torus.along.subtorus}
Let $X$ be a complete toric variety and let $T$ be a subtorus of the torus $T_X$ of $X$. 
Let $N \supseteq N_T$ be the lattices of one-parameter subgroups of $T_X$ and $T$.
Let $Y$ be a complete toric variety given by a fan in $N/N_T \otimes \mathbb{R}$ whose set of rays are the images of the rays of X under the natural projection $N \otimes \mathbb{R} \rightarrow N/N_T \otimes \mathbb{R}$ and let $e$ be a point in the torus of $Y$. 
Then, 
\[
\operatorname{Cox}(\operatorname{Bl}_{\overline{T}} X) 
\textnormal{ is finitely generated} 
\ \ \ \Longleftrightarrow \ \ \
\operatorname{Cox}\left(\operatorname{Bl}_e Y\right) 
\textnormal{ is finitely generated}. 
\]
\end{theorem}

\begin{figure}
    \centering
    \includegraphics[width=7cm]{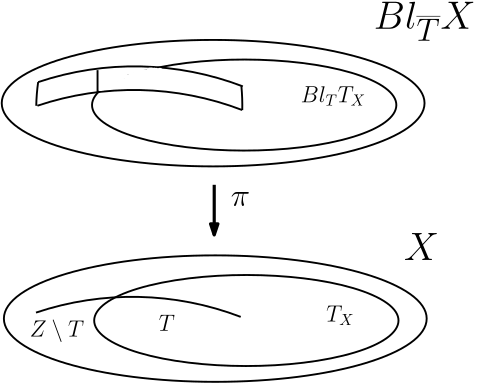}
\caption{Geometric setting in Theorem~\ref{theorem.blow-up.torus.along.subtorus}}
\label{fig:blow up}
\end{figure}

\begin{proof}
Without loss of generality we may assume that $e$ is the unit of the torus of $Y$. 
Let us denote $Z=\overline{T}$ and 
$B=\operatorname{Bl}_Z X$. 
Let 
$\pi: B \longrightarrow X$ be the blow-up morphism.
Let $U \subseteq X$ be the complement in $X$ of the $T_X$-invariant subvarieties of $X$ of codimension at least $2$.
Then $U$ is a toric variety with fan given by $\{0\}$ and the rays of $\Delta_X$.
For any algebraic set $V$ with an action of $T$ we let $\Gamma_T(V)$ denote 
the open subset of $V$ consisting of the points with zero-dimensional stabilizer.

Step 1: We show that the quotient $\Gamma_T(B)/T$ exists as a prevariety and it is equal to the union 
\begin{equation} \label{expresion.gamma.B}
\Gamma_T(\pi^{-1}(U \smallsetminus Z)) /T
\ \cup \ 
\Gamma_T(\pi^{-1}\left(T_X\right)) /T
\ \cup \ 
\Gamma_T(\pi^{-1}(Z \smallsetminus T)) /T
\ \cup \ 
\Gamma_T(\pi^{-1}(X \smallsetminus U)) /T,
\end{equation}
where 
$\Gamma_T(\pi^{-1}(U \smallsetminus Z)) /T$
and 
$\Gamma_T(\pi^{-1}\left(T_X\right)) /T$
are open in $\Gamma_T(B)/T$, and  
$\Gamma_T(\pi^{-1}(Z \smallsetminus T)) /T$
and 
$\Gamma_T(\pi^{-1}(X \smallsetminus U)) /T$
are closed in $\Gamma_T(B)/T$.

Since $T \subseteq T_X$ it is clear that $X=(U \smallsetminus Z) \cup T_X \cup (Z \smallsetminus T) \cup(X \smallsetminus U)$.
This expresses $X$ as the union of four $T$-invariant subsets, the first two open in $X$ and the latter two closed in $X$. 
By taking their inverse images under $\pi$ we get  
\begin{equation} \label{eqn.decomposition.B}
B=\pi^{-1}(U \smallsetminus Z) \ \cup \ \pi^{-1}\left(T_X\right) \ \cup \ \pi^{-1}(Z \smallsetminus T) \ \cup \ \pi^{-1}(X \smallsetminus U).
\end{equation}
The four sets on the right-hand side of (\ref{eqn.decomposition.B}) are $T$-invariant and we can obtain their subsets of points with zero-dimensional stabilizers by intersecting each one of them with $\Gamma_T(B)$, so we get 
\begin{equation} \label{eqn.decomposition.gamma.B}
\Gamma_T(B)=\Gamma_T(\pi^{-1}(U \smallsetminus Z)) 
\ \cup \ 
\Gamma_T(\pi^{-1}\left(T_X\right)) 
\ \cup \ 
\Gamma_T(\pi^{-1}(Z \smallsetminus T)) 
\ \cup \ 
\Gamma_T(\pi^{-1}(X \smallsetminus U)).
\end{equation} 
The first two sets on the right-hand side of (\ref{eqn.decomposition.gamma.B}) are open in $\Gamma_T(B)$ and the latter two are closed in $\Gamma_T(B)$. 
The geometric quotient $\Gamma_T(B)/T$ exists as a prevariety by \cite{sumihiro1974equivariant}. 
The quotients by the action of $T$ on each of the four sets on the right-hand side of (\ref{eqn.decomposition.gamma.B}) exist as unions of prevarieties as well 
since geometric quotients yield geometric quotients when restricted to open or closed invariant subsets (see \cite[Definition 1.2.3.1 and Proposition 1.2.3.9]{CoxRingsBook2015}).
Now, we pass to the geometric quotients of the sets in (\ref{eqn.decomposition.gamma.B}) by the action of $T$ to get the expression for $\Gamma_T(B)/T$ in (\ref{expresion.gamma.B}) satisfying the desired conditions.

Step 2: 
We show that the complement of the open subset 
\[
W:= \Gamma_T(\pi^{-1}(U \smallsetminus Z)) /T \cup \Gamma_T(\pi^{-1}\left(T_X\right)) /T
\] 
of $\Gamma_T(B)/T$ has codimension at least 2. 

By Step 1, it is enough to show that the closed subsets $\Gamma_T(\pi^{-1}(Z \smallsetminus T)) /T$ and $\Gamma_T(\pi^{-1}(X \smallsetminus U)) /T$ of $\Gamma_T(B)/T$ have codimension at least 2.
To see this, first notice that $\pi^{-1}(Z \smallsetminus T)$ and $\pi^{-1}(X \smallsetminus U)$ are closed subsets of $B$ with codimension at least 2, and hence $\Gamma_T(\pi^{-1}(Z \smallsetminus T))$ and $\Gamma_T(\pi^{-1}(X \smallsetminus U))$ are closed subsets of $\Gamma_T(B)$ with codimension at least 2.
Now, we note that the stabilizers of all points in $\Gamma_T(\pi^{-1}(Z \smallsetminus T))$, $\Gamma_T(\pi^{-1}(X \smallsetminus U))$ and $\Gamma_T(B)$ under the $T$-action are zero-dimensional. Then, taking the quotient by the $T$-action reduces the dimension of each set by exactly $\operatorname{dim}T$. Therefore, $\Gamma_T(\pi^{-1}(Z \smallsetminus T)) /T$ and $\Gamma_T(\pi^{-1}(X \smallsetminus U)) /T$ have codimension at least 2 in $\Gamma_T(B)/T$, and our claim in Step 2 is proven.

Step 3: We show that $W$ is isomorphic to $W'=\operatorname{Bl}_{T/T}(\Gamma_T(U)/T)$.

Let 
$W_1:=\Gamma_T(\pi^{-1}(U \smallsetminus Z)) /T$ and 
$W_2:=\Gamma_T(\pi^{-1}\left(T_X\right)) /T$, 
and hence we have an open cover $W:=W_1 \cup W_2$. 
Similarly, let 
$W'_1$ be the complement of the point $T/T$ in $\Gamma_T(U)/T$ 
and 
$W'_2:=\operatorname{Bl}_{T/T}(T_X/T)$. 
We have that $T_X=\Gamma_T(T_X)$ is an open $T$-invariant subset of $\Gamma_T(U)$
and hence $W'_2$ is an open subset of $W'$. Then, we have an open cover $W'=W'_1 \cup W'_2$. 

The $T$-stabilizers of all points in $Z \smallsetminus T$ are positive dimensional since $\operatorname{dim}(Z \smallsetminus T)<\operatorname{dim}(Z)=\operatorname{dim}(T)$, 
and hence $\Gamma_T(Z \smallsetminus T)$ is empty. 
The blow-up morphism $\pi$ gives a $T$-equivariant identification between $\pi^{-1}(U \smallsetminus Z)$ and $U \smallsetminus Z$ 
and then $\Gamma_T(\pi^{-1}(U \smallsetminus Z))$ is identified via $\pi$ with 
\[
\Gamma_T(U \smallsetminus Z) 
= \Gamma_T(U \smallsetminus Z) \cup \Gamma_T(Z \smallsetminus T) 
= \Gamma_T(U \smallsetminus T)
= \Gamma_T(U) \smallsetminus \Gamma_T(T)
= \Gamma_T(U) \smallsetminus T.
\]
Hence, $W_1=\Gamma_T(\pi^{-1}(U \smallsetminus Z)) /T$ is identified with 
$(\Gamma_T(U) \smallsetminus T) /T = \Gamma_T(U)/T \smallsetminus T/T =W'_1$.

On the other hand, we have that $\pi^{-1}\left(T_X\right)=\operatorname{Bl}_T(T_X)$ since the blow-up is a local construction.
The blow-up map $\operatorname{Bl}_T(T_X)\rightarrow T_X$ is $T$-equivariant and the $T$-stabilizers are trivial on $T_X$ and hence the $T$-stabilizers are also trivial on $\operatorname{Bl}_T(T_X)$. 
Hence
$\Gamma_T(\operatorname{Bl}_T(T_X))=\operatorname{Bl}_T(T_X)$ and therefore 
\[
W_2
=
\Gamma_T(\pi^{-1}\left(T_X\right)) /T 
=\Gamma_T(\operatorname{Bl}_T(T_X)) /T 
=\operatorname{Bl}_T(T_X) /T,  
\]
which by Lemma~\ref{lemma.blow-up.torus.along.subtorus} can be identified with $W'_2=\operatorname{Bl}_{T/T}(T_X/T)$. 

The intersection $W_1 \cap W_2 = \Gamma_T(\pi^{-1}(T_X \smallsetminus T)) /T$
can be naturally identified with $\Gamma_T(T_X \smallsetminus T)/T=(T_X \smallsetminus T)/T= T_X/T \smallsetminus T/T$,
which we recognize to be $W'_1 \cap W'_2$. 

In summary, $W$ is covered by the open subsets $W_1$ and $W_2$ 
and similarly 
$\operatorname{Bl}_{T/T}(\Gamma_T(U)/T)$
is covered by the open subsets $W'_1$ and $W'_2$.
We have identified $W_1$ with $W'_1$ and $W_2$ with $W'_2$.
Moreover, these identifications restricted to the intersections
$W_1 \cap W_2$
and 
$W'_1 \cap W'_2$ 
agree as they are both equal to the natural isomorphism between these sets that we described. 
Therefore, $W \cong W'$ and our claim in Step 3 is proven.

Step 4: We show that $\operatorname{Bl}_e Y$ is a separation of $W\cong W'$ and complete the argument. Proceeding as before, we see that $\Gamma_T(U)/T$ is an open set of $\Gamma_T(X)/T$ with complement of codimension at least 2. 
By Example~\ref{example.toric.case}, $Y$ is a separation of both $\Gamma_T(U)/T$ and $\Gamma_T(X)/T$, and 
the point $T/T$ of $\Gamma_T(X)/T$ corresponds to the unit element $e$ of the torus of $Y$ under the separation rational map. 
It follows that $\operatorname{Bl}_e Y$ is a separation of 
$W'=\operatorname{Bl}_{T/T}(\Gamma_T(U)/T)$. 
We have shown that the open subset 
$\Gamma_T(B)$ of $B=\operatorname{Bl}_{\overline{T}} X$ has a prevariety geometric quotient $\Gamma_T(B)/T$
by its $T$-action, and $\Gamma_T(B)/T$ admits $\operatorname{Bl}_e Y$ as a separation. 
Using \cite[Theorem 1.2]{hausen2010cox}, we conclude that 
$\operatorname{Cox}(\operatorname{Bl}_{\overline{T}} X)$ 
is finitely generated 
if and only if $\operatorname{Cox}\left(\operatorname{Bl}_e Y\right)$ 
is finitely generated, as desired. 
\end{proof}

\subsection{The Mori dream space property of $T_{d,n}$} \label{subsection.non-MDS.Tdn}
To conclude this section, we use the results from the previous subsection to prove Theorem~\ref{Main thm: II}. Then, we prove Theorem~\ref{theorem.Tdn.not.MDS} in Theorem~\ref{theorem.Tdn.not.MDS.chapter.4}, where we show that $T_{d,n}$ is a Mori dream space for $n\leq 3$ and that it is not a Mori dream space for $n\geq 9$. 

\begin{theorem}[Theorem~\ref{Main thm: II}] \label{thm:MainLM.Part.b}
There is an irreducible closed locus $\delta \subset T_{d,n}^{LM}$ such that the canonical map $T_{d,n} \longrightarrow T_{d,n}^{LM}$ factors as 
\begin{align*}
T_{d,n} \longrightarrow \operatorname{Bl}_{\delta} T_{d,n}^{LM}
\longrightarrow T_{d,n}^{LM}
\end{align*}
and $\operatorname{Bl}_{\delta} T_{d,n}^{LM}$ is not a Mori dream space for $n \geq 9$.
\end{theorem}

\begin{proof}
We use the notation introduced in 
Section~\ref{sec:PreliminaryToricModel}.
Let $N= \mathbb{Z}^{d(n-1)} \! / ( \sum_{i=1}^{n-1} \sum_{k=1}^{d} e^k_i =0 )$ be the lattice of one-parameter subgroups of both $\mathbb{P}^{d(n-1)-1}$ and $X=T_{d,n}^{LM}$. 
By Corollary~\ref{toricTdn}, $T_{d,n}^{LM}$ is a toric variety associated to a fan $\Delta_{T_{d,n}^{LM}}$ in $N \otimes \mathbb{R} = \mathbb{R}^{d(n-1)} \! / ( \sum_{i=1}^{n-1} \sum_{k=1}^{d} e^k_i =0 )$ whose rays are generated by the vectors 
\[
\left\{ \left. \overline{e}^{k}_{i} \ \right| \ 1 \leq i \leq n-1, 1 \leq k \leq d \right\}
\ \cup \
\left\{
\left.
\sum_{i \in I} \left( \overline{e}^1_i+ \ldots+ \overline{e}^d_i \right) 
\ \right| \
1 \leq |I| \leq n-2, \ 
I \subsetneq \{1, \ldots, n-1 \}
\right\}.
\]

Let $\mathcal{A'}:=\left(\frac{1}{n-2},\dots,\frac{1}{n-2},1\right)\in \mathcal{D}^T_{n}$ and $I=\{1,2,\ldots,n-1\}$. We notice that $T_{d,n}^{\mathcal{A'}}=\operatorname{Bl}_{\delta} T_{d,n}^{LM}$ where $\delta$ is the strict transform in $T_{d,n}^{LM}$ of the subvariety $\delta_{d,I}$ of $\mathbb{P}^{d(n-1)-1}$ defined by 
\begin{equation} \label{equations.subtorus}
  \delta_{d,I}= 
\bigcap_{i,j \in I}V((x_{i1}-x_{j1},\ldots , x_{id}-x_{jd})),
\end{equation} 
where $x_{ik}$ are the coordinates described in Lemma~\ref{lemma:PPCompactification}. The intersection of $\delta_{d,I}$ with the torus of $\mathbb{P}^{d(n-1)-1}$ is a subtorus $T$ defined by the same equations as $\delta_{d,I}$. 
Since $\mathbb{P}^{d(n-1)-1}$ and $X=T_{d,n}^{LM}$ share the same torus 
$T_{X}$, then $\delta$ is the closure in the toric variety $T_{d,n}^{LM}$ of the subtorus $T$ of the torus $T_{X}$.
Then, we are in the setting of Theorem~\ref{theorem.blow-up.torus.along.subtorus} with $X=T_{d,n}^{LM}$ and we proceed to identify the toric variety $Y$.

From the explicit equations in (\ref{equations.subtorus}) defining $T$, we see that the lattice of one-parameter subgroups $N' \subseteq N$ of $T$ is generated by 
\[
\left\{ \left.\sum_{i=1}^{n-1} \overline{e}^{\, k}_i \ \right| \ 1 \leq k \leq d \right\}
\ \subseteq \
N= \mathbb{Z}^{d(n-1)} \!\!\! \left/ \left( \sum_{i=1}^{n-1} \sum_{k=1}^{d} e^k_i =0 \right). \right.
\]
Therefore, the lattice $N/N'$ can be identified with 
\begin{eqnarray*}
  \mathbb{Z}^{d(n-1)} \! \left/ 
\left( \left. \sum_{i=1}^{n-1} e^k_i =0 \ \right| \ 1 \leq k \leq d \right). \right.
\end{eqnarray*}

For each $1 \leq i \leq n-1$ and $1 \leq k \leq d$ let $\widehat{e}^{\, \, k}_i$ be the image of $e^k_i \in \mathbb{Z}^{d(n-1)}$ in the quotient $N/N'=\mathbb{Z}^{d(n-1)} \! \left/ 
\left( \left. \sum_{i=1}^{n-1} e^k_i =0 \ \right| \ 1 \leq k \leq d \right) \right.$. 
Hence the images of the rays of $\Delta_{T_{d,n}^{LM}}$ in $N/N' \otimes \mathbb{R}$ are generated 
by the vectors 
\begin{eqnarray}  \label{equation.rays}
  \left\{ \left. \widehat{e}^{\, \, k}_i \ \right| \ 1 \leq i \leq n-1, \ 1 \leq k \leq d \right\}
\
\cup
\
\left\{
\left.
\sum_{i \in I} \left( \widehat{e}^{\, \, 1}_i+ \ldots+ \widehat{e}^{\, \, d}_i \right) 
\ \right| \
1 \leq |I| \leq n-2, \ 
I \subsetneq \{1, \ldots, n-1 \}
\right\}.
\end{eqnarray}

We recognize these rays as precisely the set of rays in the fan of the 
projective toric variety $\overline{P}^{LM}_{d,n+d}$, introduced by two of the authors in \cite{GR} as a compactification of the moduli space of $n+d$ distinct labeled points in $\mathbb{P}^d$, up to the usual action of $\SL_{d+1}$.
By Theorem~\ref{theorem.blow-up.torus.along.subtorus}, the Cox ring of $\operatorname{Bl}_{\delta} T_{d,n}^{LM}$ is finitely generated if and only if the Cox ring of the blow-up $\operatorname{Bl}_e \overline{P}^{LM}_{d,n+d}$ is finitely generated, where $e$ denotes the unit of the torus.

It was shown by three of the authors in \cite[Lemma 3.4]{Fulton-MacPherson2022} that the Cox ring of the blow-up $\operatorname{Bl}_e \overline{P}^{LM}_{d,n+d}$ is not finitely generated for $n \geq 9$. 
Therefore, the Cox ring of $\operatorname{Bl}_{\delta} T_{d,n}^{LM}$ is not finitely generated and it is not a Mori dream space for $n \geq 9$.  
\end{proof}

As an application, we study the Mori dream space property of the Chen-Gibney-Krashen compactification $T_{d,n}$. 

\begin{theorem}[Theorem~\ref{theorem.Tdn.not.MDS}] \label{theorem.Tdn.not.MDS.chapter.4}
The space $T_{d,n}$ is a Mori dream space for $n \leq 3$ and it is 
not a Mori dream space for $n \geq 9$.     
\end{theorem}
\begin{proof}
Since there exists a surjective morphism $T_{d,n} \rightarrow \operatorname{Bl}_{\delta} T_{d,n}^{LM}$ and $\operatorname{Bl}_{\delta} T_{d,n}^{LM}$ is not a MDS for $n \geq 9$ by Theorem~\ref{thm:MainLM.Part.b}, then $T_{d,n}$ is not a MDS for $n \geq 9$. We have that $T_{d,2} \cong \mathbb{P}^{d-1}$ is a MDS, so to complete the proof we show that 
$T_{d,3}$ is a MDS. Since, $T_{1,3} \cong \mathbb{P}^1$, we can  suppose that $d \geq 2$. 

We consider the construction presented in the proof of Theorem~\ref{thm:MainLM.Part.b}, but now specialized to the case $n = 3$.
Notice that in this case $T_{d,3} = \operatorname{Bl}_{\delta} T_{d,3}^{LM}$ since both are equal to the the space ${T_{d,n}^{\mathcal{A}}}$ 
from Theorem~\ref{theorem.wtdn} with respect to $\mathcal{A}=\{1,1,1\}$. 

Just as in the general case in the proof of Theorem~\ref{thm:MainLM.Part.b}, we deduce by Theorem~\ref{theorem.blow-up.torus.along.subtorus} that $T_{d,3}$ is a MDS if and only if
the blow up at the unit of the torus of any projective toric variety $Y$ whose rays are generated by the vectors in equation in (\ref{equation.rays}) is a MDS.   
Specializing (\ref{equation.rays}) to the case $n=3$ we find that $Y$ is defined by a fan in 
$\mathbb{Z}^{2d} \! \left/ 
\left( \left.  e^k_1 + e^k_2 =0 \ \right| \ 1 \leq k \leq d \right) \right.$ with rays generated by   
\begin{eqnarray*}
  \left\{ \left. \widehat{e}^{\, \, k}_i \ \right| \ 1 \leq i \leq 2, \ 1 \leq k \leq d \right\}
\
\cup
\
\left\{
 \widehat{e}^{\, \, 1}_1+ \ldots+ \widehat{e}^{\, \, d}_1, \
 \widehat{e}^{\, \, 1}_2+ \ldots+ \widehat{e}^{\, \, d}_2 
\right\}.
\end{eqnarray*}
We recognize these rays as the set of rays in the fan of the 
projective toric variety $(\mathbb{P}^1)^d$ blown up at two particular torus invariant points.  
Recall that the Mori dream space property is invariant up to small modifications of normal, projective, $\mathbb{Q}$-factorial varieties, so, for the purposes of studying the MDS property, we can consider any toric variety with the same set of rays.  
Therefore, we can assume that $Y$ is the blow up of $(\mathbb{P}^1)^d$ at two points. 
Since any triplet of ordered distinct points in $(\mathbb{P}^1)^d$ can be sent to any other such triplet via an automorphism of $(\mathbb{P}^1)^d$, we deduce that $T_{d,3}$ is a Mori dream space if and only if the blow up of of $(\mathbb{P}^1)^d$ along any three distinct points is a Mori dream space.  
But this last blow up is a Mori dream space by \cite[Theorem 1.3]{castravet2006hilbert} and this completes the proof.  
\end{proof}


\section{Higher-dimensional Losev-Manin spaces as fibrations of classical Losev-Manin spaces}
\label{sec:Fibrations}

In this section we complete the proof of Theorem~\ref{Main thm: III}. We first state and prove some preliminary lemmas. 
With the notation of Lemma~\ref{lemma:PPCompactification}, 
let us consider the rational map
\begin{equation}\label{rational}\PP^{d(n-1)-1}\dashrightarrow(\PP^{d-1})^{n-1}
\end{equation}
that sends the point 
\[
[x_{11}: x_{12}: \ldots : x_{1d} : x_{21}: x_{22} : \ldots x_{2d}:\ldots : x_{(n-1)1}: x_{(n-1)2} : \ldots : x_{(n-1)d}]\in \mathbb{P}^{d(n-1)-1}
\]
to the tuple 
\[
([x_{11}: x_{12}: \ldots : x_{1d}],[x_{21}: x_{22}: \ldots : x_{2d}],\ldots,[x_{(n-1)1}: x_{(n-1)2}: \ldots : x_{(n-1)d}]) \in (\PP^{d-1})^{n-1}. 
\]  
Observe that the indeterminacy locus of the above map is precisely the union of the varieties 
$\delta_{d,I}~=~\bigcap_{i \in I \setminus n} V(x_{i1},\ldots,x_{id})$ such that $|I|=2$ and $n\in I$ (see Equation \eqref{deltas}).   
Let $\bigcup\{\delta_{d,I}\, |\, n\in I\text{ and }|I|=2\}$ be the subscheme of $\mathbb{P}^{d(n-1)-1}$ given by the product of the ideal sheaves of $\delta_{d,I}$ in $\mathbb{P}^{d(n-1)-1}$ over the indices $I$  such that $|I|=2$ and $n\in I$.  

To state the next result we introduce some notation. Let $\pi_i:(\PP^{d-1})^{n-1}\to \PP^{d-1}$ denote the projection onto the $i$\textsuperscript{th} factor and let us introduce the following locally free sheaf and its locally free subsheaves 
\[
\mathcal{E}_{d,n}:=\bigoplus\limits_{i=1}^{n-1} \pi_i^* \mathcal{O}(1)
\quad\text{and}\quad
\mathcal{E}_{d,n}^{(j)}:=\bigoplus_{\substack{i=1\\i\neq j}}^{n-1} \pi_i^* \mathcal{O}(1), 
\]
for $1 \leq j \leq n-1$. Further, recall the following definition. 

\begin{definition} Let $X$ be a scheme, let $Y$ and $Z$ be subschemes of $X$, and let $\phi:\Bl_YX\to X$ be the blowup of $X$ along $Y$. The total transform of $Z$ in $\Bl_YX$ is the scheme theoretic inverse $\phi^{-1}(Z)$.
\end{definition}

We may now state the following lemma.  
If $\mathcal{E}$ is a vector bundle of over $X$, then we will denote the associated projective space bundle by $\mathbb{P}_X \left( \mathcal{E} \right)$.  
\begin{lemma}\label{firststep} 

\begin{enumerate} 
\item The blow-up of $\mathbb{P}^{d(n-1)-1}$ along $\bigcup\{\delta_{d,I}\, |\, n\in I\text{ and }|I|=2\}$ is isomorphic to the projective bundle 
\[
    \PP_{(\PP^{d-1})^{n-1}}(\mathcal{E}_{d,n}).
\]
\item Under the above isomorphism, for each $i \in \{1,\dots,n-1\}$, the total transform of $\delta_{d,\{i,n\}}$ is identified with the projective sub-bundle 
\[
    \PP_{(\PP^{d-1})^{n-1}}(\mathcal{E}_{d,n}^{(i)}).
\]
\end{enumerate}
\end{lemma}
In proving the above lemma, we use the following result which is probably well-known, yet we provide a proof for the reader's convenience.
 
\begin{lemma} \label{aux} Let $V$ be a ($k+1$)-dimensional vector subspace of a given ($m+1$)-dimensional vector space $W$. The blow-up of the projective space $\PP(W)$ along $\PP(V)$ is identified with the projective bundle \[\PP_{\PP(W/V)}(\mathcal{O}^{k+1}\oplus \mathcal{O}(1)).\]
 \end{lemma}
 \begin{proof}
 We may choose projective coordinates 
 \[
    [x_0:x_1:\ldots:x_m]
\]
 for $\PP^m=\PP(W)$ so that the subspace $\PP^k=\PP(V)$ is given by the vanishing of $x_{k+1},\ldots, x_m$. Then, the blow-up of $\PP(W)$ along $\PP(V)$ is identified with the subvariety of $\PP^m\times \PP^{m-k-1}$ defined by the vanishing of the ideal generated by $x_iy_j-x_jy_i,\,i,j=k+1,\ldots,m$, where $y_{k+1},\ldots,y_m$ are the coordinates of $\PP^{m-k-1}$. 
 
 At the same time, we may embed the bundle $\PP_{\PP^{m-k-1}}(\mathcal{O}^{k+1}\oplus\cO(1))=\Proj(\Sym(\cO^{k+1}\oplus \cO(1)))$ in $\PP^m\times\PP^{m-k-1}=\Proj(\cO_{\PP^{m-k-1}}[x_0,\ldots,x_m])$ as follows. We define a morphism of sheaves of $\cO_{\PP^{m-k-1}}$-algebras 
 $$\cO_{\PP^{m-k-1}}[x_0,\ldots,x_m]\to \Sym(\cO^{k+1}\oplus \mathcal{O}(1))$$ by the rule
 $$x_i\mapsto
 \begin{cases} e_i,\,\,\,\,i=0,\ldots k\\
 y_i,\,\,\,\,i=k+1,\ldots, m
 \end{cases}$$
 where $e_0,e_1,\dots,e_k$ is the canonical basis of $H^0(\PP^{m-k-1},\mathcal{O}_{\PP^{m-k-1}}^{k+1})$.
 This is clearly an epimorphism with kernel equal to the ideal generated by $x_iy_j-x_jy_i,\,i,j=k+1,\ldots,m$, hence we arrive at the required identification. 
 \end{proof}

\begin{proof}[Proof of Lemma~\ref{firststep}]
 One can check that the set $\{\delta_{d,I}\, |\, n\in I\text{ and }|I|=2\}$ of subvarieties of $\mathbb{P}^{d(n-1)-1}$ is a building set in the sense of Li \cite{LiLi}.
 Therefore, the blow-up of $\mathbb{P}^{d(n-1)-1}$ along $\bigcup\{\delta_{d,I}\, |\, n\in I\text{ and } |I|=2\}$ is the wonderful compactification of the above building set \cite{LiLi}, namely, it is the closure of the natural locally closed embedding
 \begin{equation}\label{closure}
  \mathbb{P}^{d(n-1)-1}\setminus\bigcup_{\substack{|I|=2\\ n\in I}} \delta_{d,I}\xhookrightarrow{}\prod_{\substack{|I|=2\\ n\in I}} \Bl_{\delta_{d,I}}\PP^{d(n-1)-1}.
 \end{equation}
 
 By Lemma~\ref{aux}, for each $i=1,\dots,n-1$, we have an isomorphism  
 $$\theta_i : \PP_{ \PP^{d-1}}({\mathcal{O}}^{d(n-2)}\oplus \mathcal{O}(1)) \rightarrow \Bl_{\delta_{d,\{i,n\}}}\PP^{d(n-1)-1},$$ 
and we denote the corresponding blowup morphism by $\beta_i:\Bl_{\delta_{d,\{i,n\}}}\PP^{d(n-1)-1} \rightarrow \PP^{d(n-1)-1}$.
  
 We may define a surjective morphism of sheaves $\phi_i:{\mathcal{O}}^{d(n-2)}\oplus \pi_i^*\mathcal{O}(1)\to \mathcal{E}_{d,n} $ as follows.
 Let $x_1,\dots x_d$ be a basis for the space of global sections of $\mathcal{O}_{\PP^{d-1}}(1)$ and $x_{i1},\dots x_{id}$ be the induced basis of global sections for the sheaf $\pi_i^*\mathcal{O}_{\PP^{d-1}}(1)$.
We define 
 $$\xi: \mathcal{O}_{\PP^{d-1}}^{d}\to \mathcal{O}_{\PP^{d-1}}(1)$$
 via $e_i\mapsto x_i$
 where $e_i$ is the $i$\textsuperscript{th} standard element of the basis of $\cO_{\PP^{d-1}}^{d}$. Then, the morphism induced by pulling back via $\pi_i$
  \[
    \pi_i^*\xi: \pi_i^*\mathcal{O}_{\PP^{d-1}}^{d}\to \pi_i^* \mathcal{O}_{\PP^{d-1}}(1)
  \]
sends $\pi_i^* e_j$ to $x_{ij}$.

We can then define the map $\phi_i$ via the following commutative diagram:
\[
    \begin{tikzcd}[column sep=6em]
        {\mathcal{O}}_{(\PP^{d-1})^{n-1}}^{d(n-2)}\oplus\pi_i^*\mathcal{O}_{\PP^{d-1}}(1) \ar[r,"\phi_i"]\ar[d,equal] & 
        \mathcal{E}_{d,n}
         \ar[d,equal] \\
        \bigoplus\limits_{j\neq i}\pi_j^*\mathcal{O}_{\PP^{d-1}}^{d} \oplus \pi_i^* \mathcal{O}_{\PP^{d-1}}(1) \ar[r, swap, "((\pi_j^*\xi)_{j\neq i}{,}\,\,id)"] &
        \bigoplus\limits_{j\neq i} \pi_j^* \mathcal{O}_{\PP^{d-1}}(1) \oplus \pi_i^* \mathcal{O}_{\PP^{d-1}}(1)
    \end{tikzcd}
\]
The above morphism $\phi_i$ is clearly surjective and thus gives rise to a closed immersion
\begin{equation}\label{immersion}
\Phi_i:\PP_{(\PP^{d-1})^{n-1}}(\mathcal{E}_{d,n})\hookrightarrow \PP_{(\PP^{d-1})^{n-1}}({\mathcal{O}}^{d(n-2)}\oplus \pi_i^*\mathcal{O}(1))
\end{equation}
which, composed with the projection (i.e., the change of basis pullback of bundles over $\pi_i:(\PP^{d-1})^{n-1}\to \PP^{d-1}$), 
$$\PP_{(\PP^{d-1})^{n-1}}({\mathcal{O}}^{d(n-2)}\oplus \pi_i^*\mathcal{O}(1))\to \PP_{\PP^{d-1}}({\mathcal{O}}^{d(n-2)}\oplus \mathcal{O}(1))$$
gives rise to a morphism 
$$\Psi_i:\PP_{(\PP^{d-1})^{n-1}}(\mathcal{E}_{d,n})\to\PP_{\PP^{d-1}}({\mathcal{O}}^{d(n-2)}\oplus \mathcal{O}(1)). $$



For each $k=1,2,\ldots,n-1$, let us define 
\[
c_k:=\beta_k \circ \theta_k \circ \Psi_k :  
\mathbb{P}_{(\mathbb{P}^{d-1})^{n-1}}(\mathcal{E}_{d,n})
\rightarrow
\PP^{d(n-1)-1}. 
\] 
Next, we fix an index $k \in \{1,2,\ldots,n-1\}$ and we will describe $c_k$ in local coordinates.

\medskip

\noindent\textbf{Choosing local coordinates.}\; 
We recall that, following the convention of \cite[Section II.7]{Hartshorne}, a point in the projectivization $\mathbb{P}(\mathcal{E})$ of a vector bundle $\mathcal{E}$ corresponds to a point $x$ on the base together with a one-dimensional quotient of the fiber $\mathcal{E}_x$ 
(that is, a subspace  $V \subseteq \mathcal{E}_x$ such that $\operatorname{dim}(\mathcal{E}_x/V)=1$).

Each point $Q \in \mathbb{P}_{(\mathbb{P}^{d-1})^{n-1}}(\mathcal{E}_{d,n})$ projects to a point $P=(P_1,\ldots,P_{n-1})\in(\mathbb{P}^{d-1})^{n-1}$ on the base. 
Given $i \in \{1,\ldots,n-1\}$,
if we write $P_i=[x_{i1}:\ldots:x_{id}] \in \mathbb{P}^{d-1}$, then there exists $j_i \in \{1,\ldots,d\}$ such that $x_{ij_i} \neq 0$.
By relabeling coordinates of each factor $\mathbb{P}^{d-1}$ if necessary, we will assume that $j_i=d$ for each $i$.  
Consequently, we fix the affine chart
\[
  U \;=\; \prod_{i=1}^{n-1}\{x_{id}\neq0\}\;\subseteq\;(\mathbb{P}^{d-1})^{n-1},
  \qquad 
  u_{ij}\,:=\,\frac{x_{ij}}{x_{id}}\quad(1\le i\le n-1,\ 1\le j\le d-1).
\]
For notational convenience we define $u_{id}:=1$ for each $1\leq i \leq n-1$, and every point of $U$ can be written as
\[
  P=(P_{1},\dots,P_{n-1}),
  \qquad
  P_{i}=[u_{i1}:\dots:u_{i(d-1)}:1]\in\mathbb{P}^{d-1}.
\]
Over $U$, each summand $\pi_{i}^{*}\mathcal{O}_{\mathbb{P}^{d-1}}(1)$
of the bundle
$
  \mathcal{E}_{d,n}=\bigoplus_{i=1}^{n-1}\pi_{i}^{*}\mathcal{O}(1)
$
is trivialized by the section $x_{id}$, so that
$
  \mathcal{E}_{d,n}|_{U}\cong\mathcal{O}_{U}^{\,n-1}.
$
This isomorphism gives us a basis $(e_{1},\dots,e_{n-1})$ for that free module, and in turn gives us a basis for the fiber of $(\mathcal{E}_{d-n})_P$ over each point of $P$ of $U$, which for simplicity will also be denoted as $(e_{1},\dots,e_{n-1})$. 
The point $Q$ is determined by its projection $P$ to the base and 
a one-dimensional quotient
$
  (\mathcal{E}_{d,n})_P \twoheadrightarrow \mathbb{C}
$
which, up to scaling, corresponds to the functional
$
  \ell=a_{1}e_{1}^{\vee}+\dots+a_{n-1}e_{\,n-1}^{\vee},
$
i.e., by the homogeneous coordinates
$
  [a_{1}:\dots:a_{n-1}]\in\mathbb{P}^{\,n-2}.
$
From here on, we will call $(P_1,\ldots,P_{n-1},[a_{1}:\dots:a_{n-1}])$ the local coordinates of $Q$. 
 
\medskip
\noindent\textbf{The morphisms $\Phi_{k}$ and $\Psi_{k}$ in local coordinates.}\;
In our local coordinates, the surjection
\[
  \mathcal{O}^{\,d(n-2)}\!\oplus\!\pi_{k}^{*}\mathcal{O}(1)
  \;\twoheadrightarrow\;
  \mathcal{E}_{d,n}
\]
acts on the fibers as follows.  
Write a point of the domain fiber over $P$ as
\[
  (\,v_{1},\dots,v_{k-1},v_{k+1},\dots,v_{n-1},\lambda\,) \in \mathbb{C}^{d(n-2)+1},
  \quad v_{i}=(v_{i1},\dots,v_{id})\in \mathbb{C}^{d},\quad 
        \lambda\in \mathbb{C},
\]
omitting the $i=k$ block.  The image of this point in
$(\mathcal{E}_{d,n})_P \cong \mathbb{C}^{\,n-1}$ is
\[
\begin{array}{l}
\left( u_{11}v_{11} + \dots + u_{1d}v_{1d},\;
\dots,\;
u_{(k-1)1}v_{(k-1)1} + \dots + u_{(k-1)d}v_{(k-1)d},\;
\lambda,\;  \right. \\[2pt]
\hspace{4em} \left. 
u_{(k+1)1}v_{(k+1)1} + \dots + u_{(k+1)d}v_{(k+1)d},\; 
 \dots,\;
u_{(n-1)1}v_{(n-1)1} + \dots + u_{(n-1)d}v_{(n-1)d} \right).
\end{array}
\]
We can compute the fiber coordinates of $\Phi_{k}\!\bigl(
    P,[a_{1}:\dots:a_{\,n-1}]\bigr)$ by evaluating the functional $\ell$ on this image.   
It follows that on the chart $U$, $\Phi_{k}\!\bigl(
    P,[a_{1}:\dots:a_{\,n-1}]
  \bigr)$ is given by 
\[
  \bigl(
    P,\,
    [\,a_{1}u_{11}:\dots:a_{1}u_{1d}:\;\dots\;:
       a_{k-1}u_{\,(k-1)d}:\,a_{k}:\,a_{\,k+1}u_{\,(k+1)1}:\dots:
       a_{\,n-1}u_{\,(n-1)d}\,]
  \bigr)
\]
where the bracketed list $[\cdots]$ is a point of
$\mathbb{P}^{\,d(n-2)}$ lying in the fiber of
$\mathbb{P}_{(\mathbb{P}^{d-1})^{n-1}}\bigl(\mathcal{O}^{d(n-2)}\oplus \pi^{*}_{k}\mathcal{O}(1)\bigr)$
over $P$.
Therefore, on the chart $U$, $\Psi_{k}\!\bigl(
    P,[a_{1}:\dots:a_{\,n-1}]
  \bigr)$ is given by 
\[
  \bigl(
    P_{k},\,
    [\,a_{1}u_{11}:\dots:a_{1}u_{1d}:\;\dots\;:
       a_{k-1}u_{\,(k-1)d}:\,a_{k}:\,a_{\,k+1}u_{\,(k+1)1}:\dots:
       a_{\,n-1}u_{\,(n-1)d}\,]
  \bigr)
\]
where the bracketed list $[\cdots]$ is now a point of
$\mathbb{P}^{\,d(n-2)}$ lying in the fiber of
$\mathbb{P}_{\mathbb{P}^{d-1}}\bigl(\mathcal{O}^{d(n-2)}\oplus\mathcal{O}(1)\bigr)$
over $P_{k}$.

\medskip
\noindent\textbf{From \(\Psi_{k}\) to \(c_{k}\) via the incidence model of the blow-up.}

\smallskip

\noindent

Write the coordinates on \(\mathbb{P}^{d(n-1)-1}\) in \(d\)-tuples
\[
[x_{11}:\dots:x_{1d}:\;\dots\;:\;x_{\,(n-1)1}:\dots:x_{\,(n-1)d}],
\]
so the center
\(
\delta_{d,\{k,n\}}\simeq\mathbb{P}^{d-1}
\)
is defined by \(x_{k1}=\dots=x_{kd}=0\). 
Using the incidence description,   
Lemma~\ref{aux} identifies the blow-up of $\mathbb{P}^{d(n-1)-1}$ along $\delta_{d,\{k,n\}}$ with
\begin{equation}\label{equation.bilinear.relations}
\operatorname{Bl}_{\delta_{d,\{k,n\}}}\mathbb{P}^{d(n-1)-1}
=\Bigl\{(x,y)\in\mathbb{P}^{d(n-1)-1}\times\mathbb{P}^{d-1}\;\Bigm|\;
       x_{k\ell}\,y_{k\ell'}=x_{k\ell'}\,y_{k\ell}\;\;
       \forall\,\ell,\ell'\Bigr\}, 
\end{equation}
where $y_{k1},\ldots,y_{kd}$ are the coordinates of $\mathbb{P}^{d-1}$ and the blow-down is the projection
\(
\beta_{k}(x,y)=x.
\)

\smallskip
Let us describe the isomorphism \(\theta_{k}\) in Lemma~\ref{aux} locally. 
Let
\(P_{k}=[u_{k1}:\dots:u_{k(d-1)}:1]\in\mathbb{P}^{d-1}\)
be the \(k\)-th factor of the base point \(P\in U\).
On the chart \(U\) a point of the fibre of  
\(\mathbb{P}_{\mathbb{P}^{d-1}}
   (\mathcal{O}^{d(n-2)}\oplus\mathcal{O}(1))\)
is
\[
[b_{11}:\dots:b_{1d}:\;\dots\;:b_{\,(k-1)d}:\lambda:
  b_{\,(k+1)1}:\dots:b_{\,(n-1)d}].
\]
Define
\[
\theta_{k}\bigl(P_{k},[b_{\bullet},\lambda]\bigr)=
\bigl(x,y\bigr)\quad\text{with}\quad
\left\{
\begin{aligned}
x_{k\ell}&=\lambda\,u_{k\ell}\;&(\ell<d),\\
x_{kd}&=\lambda,\\
x_{j\ell}&=b_{j\ell}\;& (j\neq k),\\
y_{k\ell}&=u_{k\ell}\;&(\ell<d),\qquad
y_{kd}=1.
\end{aligned}
\right.
\]
These \(x\) and \(y\) satisfy the bilinear relations of \eqref{equation.bilinear.relations},
so \(\theta_{k}\) lands in the blow-up. 

\smallskip

We have seen that the point $\Psi_{k}\bigl(P,[a_{1}:\dots:a_{n-1}]\bigr)$ is given by
\[
  \bigl(
    P_{k},\,
    [\,a_{1}u_{11}:\dots:a_{1}u_{1d}:\;\dots\;:
       a_{k-1}u_{\,(k-1)d}:\,a_{k}:\,a_{\,k+1}u_{\,(k+1)1}:\dots:
       a_{\,n-1}u_{\,(n-1)d}\,]
  \bigr)
\]
so \(\lambda=a_{k}\).
Applying \(\theta_{k}\) replaces the entry \(a_{k}\) by the block
\([a_{k}x_{k1}:\dots:a_{k}x_{kd}]\), and  
applying \(\beta_{k}\) then drops \(y\) and keeps the \(x\)-coordinates,
giving

\begin{align} \label{equation.ck}
\left[\, a_{1}u_{11}:\dots:a_{1}u_{1d}:\dots :
  a_{k}u_{k1}:\dots:a_{k}u_{kd}:\dots:\;
  a_{\,n-1}u_{\,(n-1)1}:\dots:a_{\,n-1}u_{\,(n-1)d}\,\right]
\in \mathbb{P}^{d(n-1)-1}. 
\end{align} 

\medskip
\noindent\textbf{Conclusion about $c_k$.}  
Given a point $Q$ in $\mathbb{P}_{(\mathbb{P}^{d-1})^{n-1}}$ mapping to the base to the point
$P=(P_1,\ldots,P_{n-1})$
with 
$P_{i}=[u_{i1}:\dots:u_{i(d-1)}:1]\in\mathbb{P}^{d-1}$, the answer in \eqref{equation.ck} is 
\begin{equation}\label{equation.ck2}
c_{k}(Q) = 
\left[\, a_{1}u_{11}:\dots:a_{1}u_{1d}:\dots :
  a_{k}u_{k1}:\dots:a_{k}u_{kd}:\dots:\;
  a_{\,n-1}u_{\,(n-1)1}:\dots:a_{\,n-1}u_{\,(n-1)d}\,\right]
\in \mathbb{P}^{d(n-1)-1},
\end{equation}
which provides the description of $c_k$ in local coordinates that we were looking for.


\medskip

We notice that $c_k: 
\mathbb{P}_{(\mathbb{P}^{d-1})^{n-1}}(\mathcal{E}_{d,n})
\rightarrow
\mathbb{P}^{d(n-1)-1}$ is independent of $k$, and we will denote this common morphism as $c: 
\mathbb{P}_{(\mathbb{P}^{d-1})^{n-1}}(\mathcal{E}_{d,n})
\rightarrow 
\mathbb{P}^{d(n-1)-1}$. 
Let us consider the open subset $W$ of 
\[
\mathbb{P}_{(\mathbb{P}^{d-1})^{n-1}}(\mathcal{E}_{d,n})=\mathbb{P}_{(\mathbb{P}^{d-1})^{n-1}}\left(\bigoplus\limits_{i=1}^{n-1} \pi_i^* \mathcal{O}(1) \right)
\]
consisting of the points $Q$ lying over $U \;=\; \prod_{i=1}^{n-1}\{x_{id}\neq0\}$ that expressed in local coordinates $Q=(P,[a_1:\cdots:a_{n-1}])$ satisfy that $a_i\neq 0$ for all $i$. 
We see that $W$ is well-defined since it is equal to the complement in the preimage of $U$ of the union of the prime divisors $\PP_{(\PP^{d-1})^{n-1}}(\mathcal{E}_{d,n}^{(j)}) \subseteq \PP_{(\PP^{d-1})^{n-1}}(\mathcal{E}_{d,n})$ for $j=1,\ldots,n-1$. 
Also, we define $V:=\mathbb{P}^{d(n-1)-1} \setminus \bigcup\{\delta_{d,I}\, |\, n\in I\text{ and }|I|=2\}$ of $\mathbb{P}^{d(n-1)-1}$. Notice that $c$ maps $W$ isomorphically onto the open subset 
\[
V\setminus \{ x_{1d}=0,\dots,x_{(n-1)d}=0\},
\]
 which is dense in $V$. It follows that the morphism $(\prod_{i=1}^{n-1} \Psi_i)\circ \Delta$, where $\Delta $ is the ($n-1$)-fold diagonal of $\PP_{(\PP^{d-1})^{n-1}}(\mathcal{E}_{d,n})$, maps $W \subseteq \PP_{(\PP^{d-1})^{n-1}}(\mathcal{E}_{d,n}) $ into 
\[
V = \mathbb{P}^{d(n-1)-1}\setminus\bigcup_{\substack{|I|=2\\ n\in I}} \delta_{d,I}\subseteq  \prod_{\substack{|I|=2\\ n\in I}} \Bl_{\delta_{d,I}}\PP^{d(n-1)-1},
\]
and has a dense image.  

Since $W$ is dense in $\PP_{(\PP^{d-1})^{n-1}}(\mathcal{E}_{d,n})$, we deduce that $(\prod_{i=1}^{n-1} \Psi_i)\circ \Delta$ maps $\PP_{(\PP^{d-1})^{n-1}}(\mathcal{E}_{d,n})$ to the closure $\overline{V}$ of $V$ in $\prod_{{|I|=2, n\in I}} \Bl_{\delta_{d,I}}\PP^{d(n-1)-1}$ which is the blow up of $\PP^{d(n-1)-1}$ along $\bigcup\{\delta_{d,I}\, |\, n\in I\text{ and }|I|=2\}$, as seen in \eqref{closure}. In particular,  $\overline{V}$ is smooth.  
Therefore, to conclude the proof of part (1) of the lemma, it suffices to show that the morphism 
\begin{equation}
    (\prod_{i=1}^{n-1} \Psi_i)\circ \Delta: \PP_{(\PP^{d-1})^{n-1}}(\mathcal{E}_{d,n}) \rightarrow \overline{V} \label{eqn.isomorphism}
\end{equation}
is an isomorphism. 
We use the following version of Zariski's Main Theorem \cite[Theorem~12.83]{GortzWedhorn}: If $f:X\to Y$ is a bijective proper morphism of varieties, with $Y$ normal, then $f$ is an isomorphism.

We now verify that the hypotheses of this theorem hold with $X=\PP_{(\PP^{d-1})^{n-1}}(\mathcal{E}_{d,n})$, $Y=\overline{V}$ and $f=(\prod_{i=1}^{n-1} \Psi_i)\circ \Delta$. 
The variety $Y$ is smooth, in particular it is normal. 
The varieties $X$ and $Y$ are projective, in particular they are proper. 
It follows that $f$ is a proper morphism, and since it is dominant it must be surjective. 
We are reduced to show that $f$ is injective, and for this we assume that $Q,Q' \in X$ are such that $f(Q)=f(Q')$. 
Let $P=(P_1,\ldots,P_{n-1})$ and $P'=(P'_1,\ldots,P'_{n-1})$ be the respective projections of $Q$ and $Q'$ to the base $(\mathbb{P}^d)^{n-1}$. 
Since $\Psi_i(Q)=\Psi_i(Q)$ for each $i=1,\ldots,n-1$, we can further compose with the projections to each factor $\mathbb{P}^{d-1}$ and deduce that $P_i=P'_i$ for each $i=1,\ldots,n-1$. 
Therefore, $P=P'$ and hence the points $Q$ and $Q'$ are on the same fiber of the morphism $\PP_{(\PP^{d-1})^{n-1}}(\mathcal{E}_{d,n}) \rightarrow (\PP^{d-1})^{n-1}$.   
By relabeling the coordinates, if necessary, we may assume that $P,P' \in U$ and use the local coordinates that we introduced before. 

Let us say that the local coordinates of $Q$ and $Q'$ are respectively
$((P_1,\ldots,P_{n-1}),[a_1:\cdots:a_{n-1}])$ and $((P_1,\ldots,P_{n-1}),[a'_1:\cdots:a'_{n-1}])$
with 
$P_{i}=[u_{i1}:\dots:u_{i(d-1)}:1]\in\mathbb{P}^{d-1}$
and $[a_1:\cdots:a_{n-1}],[a'_1:\cdots:a'_{n-1}] \in \mathbb{P}^{n-2}$. 
For convenience we define  $u_{id}=1$ for all $i$. 
We have that $c(Q)=c(Q')$ and, from Equation~\eqref{equation.ck2}, that
\begin{align*}
c(Q) &= 
\left[\, a_{1}u_{11}:\dots:a_{1}u_{1d}:\dots :
  a_{k}u_{k1}:\dots:a_{k}u_{kd}:\dots:\;
  a_{\,n-1}u_{\,(n-1)1}:\dots:a_{\,n-1}u_{\,(n-1)d}\,\right]
\in \mathbb{P}^{d(n-1)-1}, \\
c(Q') &= 
\left[\, a'_{1}u_{11}:\dots:a'_{1}u_{1d}:\dots :
  a'_{k}u_{k1}:\dots:a'_{k}u_{kd}:\dots:\;
  a'_{\,n-1}u_{\,(n-1)1}:\dots:a'_{\,n-1}u_{\,(n-1)d}\,\right]
\in \mathbb{P}^{d(n-1)-1}.
\end{align*}
It follows that $[a_1:\cdots:a_{n-1}] = [a'_1:\cdots:a'_{n-1}]$, and hence $Q=Q'$, from which we deduce $f$ is injective. 

By Zariski's Main Theorem the morphism $f= (\prod_{i=1}^{n-1} \Psi_i)\circ \Delta$ in \eqref{eqn.isomorphism} is an isomorphism, concluding the proof of part (1) of the lemma.

\medskip


For the second part, let  $F_{d,n}^{(i)}$ be the exceptional divisor of $\Bl_{\delta_{d,\{i,n\}}}\PP^{d(n-1)-1}$  for each $i=1,\ldots,{n-1}$. Under the identification of $\Bl_{\delta_{d,\{i,n\}}}\PP^{d(n-1)-1}$ with $\PP_{ \PP^{d-1}}({\mathcal{O}}^{d(n-2)}\oplus \mathcal{O}(1))$, the divisor $F_{d,n}^{(i)}$ corresponds to the sub-bundle
$\PP_{ \PP^{d-1}}({\mathcal{O}}^{d(n-2)})$ of $\PP_{ \PP^{d-1}}({\mathcal{O}}^{d(n-2)}\oplus \mathcal{O}(1))$ induced by the projection 
${\mathcal{O}}^{d(n-2)}\oplus \mathcal{O}(1)\to {\mathcal{O}}^{d(n-2)}$. Hence, the restriction of $\PP_{ \PP^{d-1}}({\mathcal{O}}^{d(n-2)})$ under (\ref{immersion}) is precisely the sub-bundle $\PP_{(\PP^{d-1})^{n-1}}(\mathcal{E}_{d,n}^{(i)})$ of $\PP_{(\PP^{d-1})^{n-1}}(\mathcal{E}_{d,n})$. Consequently the restrictions of the divisors $$F_{d,n}^{(j)}\times \prod_{\substack{i= 1\\i\neq j}}^{n-1}\Bl_{\delta_{d,\{i,n\}}}\PP^{d(n-1)-1}$$ under the immersion
$$(\prod_{i=1}^{n-1} \Psi_i)\circ \Delta:\PP_{(\PP^{d-1})^{n-1}}(\mathcal{E}_{d,n})\to\prod_{i=1}^{n-1}\PP_{\PP^{d-1}}({\mathcal{O}}^{d(n-2)}\oplus \mathcal{O}(1))=\prod_{i=1}^{n-1}\Bl_{\delta_{d,\{i,n\}}}\PP^{d(n-1)-1} $$
coincide with the sub-bundles $\PP_{(\PP^{d-1})^{n-1}}(\mathcal{E}_{d,n}^{(j)})$ for $j=1,\dots,n-1$ and their union is precisely the boundary of the closure of (\ref{closure}). Part (2) of the lemma now follows immediately.
\end{proof}

We are now able to complete the proof of Theorem~\ref{Main thm: III}. Namely, that $T_{d,n}^{LM}$ is a Zariski locally trivial fibration over $(\PP^{d-1})^{n-1}$, with fiber isomorphic to the Losev-Manin space $\overline{M}_{0,n+1}^{LM}$.

\begin{proof}[Proof of Theorem~\ref{Main thm: III}]
By Theorem~\ref{theorem.wtdn I} and \cite[Lemma 3.2]{LiLi}, the space $T_{d,n}^{LM}$ can be obtained as the iterated blow-up of $\PP^{d(n-1)-1}$ along the \emph{total} transforms of the $\delta_{d,I}$, where $n\in I$, in any order. We may therefore obtain $T_{d,n}^{LM}$ as a sequence of blow-ups performed in the following particular order:
\begin{enumerate}
  \item blow up $\PP^{d(n-1)-1}$ successively along the total transforms of $\delta_{d,I}$, where $n\in I$ and $|I|=2$ in any order; denote the resulting blow-up by $\textbf{P}^{[1]}$ .
  \item blow up the total transforms of rest of the $\delta_{d,I}$ in $\textbf{P}^{[1]}$ (i.e., those with $|I|>2$ and $n\in I$) in any order.
\end{enumerate}

By Lemma~\ref{aux}, the blow-up of the first stage is identified with the bundle $\PP(\mathcal{O}(1,0,\dots,0)\oplus \dots \oplus\mathcal{O}(0,0,\dots,1))$ over $(\PP^{d-1})^{n-1}$. Let $U$ be an open subset of $(\PP^{d-1})^{n-1}$ where that bundle is trivial. By Lemma~\ref{aux} again, 
over $U$, the blow-up of the first stage is identified with $U\times \PP^{n-2}$ and its exceptional divisors coincide the product of $U$ with the toric boundary divisors of 
$\PP^{n-2}$. Further, note that for each 
$I$ containing $n$ and $|I|\geq 2$ we have 
\[
    \delta_{d,I}=
    \bigcap_{i \in I \setminus \{n\}}
    \delta_{d,\{i,n\}}.
\]
Consequently, the set of blown-up centers in stage (2) is identified with the set consisting of the product of $U$ with all possible intersections of the toric boundary divisors of $\PP^{n-2}$, that is the product of $U$ with the set $\{\delta_{1,I}|\,n\in I\, \text{and}\, |I|\geq 2 \}$. Therefore, the resulting blow-up is identified with the product $U\times T_{1,n}^{LM}$ by Theorem~\ref{theorem.wtdn I} and $T_{1,n}^{LM}$ is identified with the standard Losev-Manin space of rational $(n+1)$-pointed curves (Corollary \ref{toricTdn}). 
\end{proof}

Finally, let us describe an application of the results in this section. Recall that an important problem in toric geometry is to characterize the toric varieties with a full strong exceptional collection of line bundles; for the definition and further details see \cite[Definitions 1.57 and 8.31]{huybrechts2006fourier}.  
By results of Orlov, smooth varieties constructed by an iterative blow-up of the projective space along smooth centers have a full strong exceptional collection; see \cite[Section 3]{castravet2020derived}. Examples include the spaces $\overline{M}_{0,n}$ and $T_{d,n}^{LM}$. Yet, finding the values of $d$ and $n$ for which such a collection consisting of line bundles exists is most interesting. Indeed, it was believed that all toric varieties carry such collection---a statement known as \textit{King's conjecture}. 
However, counterexamples have been constructed; see \cite{hille2006counterexample}.
Next, we provide a new family of toric varieties satisfying this property.
\begin{corollary}\label{cor:DC}
    If $d \geq 1$ and $n \leq 4$, then $T_{d,n}^{LM}$ has a full strong exceptional collection of line bundles. 
\end{corollary}
\begin{proof}
By \cite[Theorem 1.3]{costa2011derived}, if $f:X \rightarrow Z$ is a Zariski locally trivial fibration of smooth, complex projective varieties, with fiber $F$ such that both $F$ and $Z$ have a full strongly exceptional collection of line bundles, then there exists a full strongly exceptional collection of line bundles on $X$.
In our case, the base of the fibration, $(\mathbb{P}^{d-1})^{n-1}$, clearly has such a collection. Therefore, the result follows from a similar one for the fiber, namely the Losev-Manin space $\overline{M}_{0,n+1}^{LM}$. 
First, we suppose that $d \geq 2$. The case $T_{d,2}^{LM} \cong \mathbb{P}^{d-1}$ is trivial. For $n=3$, it follows from $\overline{M}_{0,4}^{LM} \cong \mathbb{P}^1$ and, for $n=4$, it follows from the fact that $\overline{M}_{0,5}^{LM}$ is the smooth toric del Pezzo surface $\Bl_{3pt}\mathbb{P}^2$; see \cite[Theorem~7.3]{borisov2009conjecture} and \cite{king1997tilting}. For the case $d=1$, it follows from the isomorphisms $T_{1,4}^{LM} \cong \overline{M}_{0,5}^{LM}$ and $T_{1,3}^{LM} \cong \overline{M}_{0,4}^{LM}$.
\end{proof}


\vspace{5mm}

\bibliographystyle{amsalpha}
\bibliography{biblioTopicsTdn} 

\vspace{10mm}

\end{document}